\documentclass[twoside,11pt]{article}

\newcommand{\techrep}{}

\usepackage{jmlr2e}
\usepackage{ifthen}
\usepackage{amsmath}
\usepackage{url}
\usepackage{subfigure}
\usepackage{pgf}
\usepackage{tikz}
\usetikzlibrary{arrows}

\newcommand{\ie}{i.e., }
\newcommand{\eg}{e.g., }

\newcommand{\qedsymbol}{\BlackBox}
\newcommand{\qed}{\hfill\qedsymbol}
\newcommand{\qedtag}[1]{\tag*{\lower#1\hbox{\qedsymbol}}}
\newenvironment{proofsketch}{\noindent \textbf{Proof sketch\hskip.5em\relax}}{\qed\vskip\baselineskip}
\newtheorem{theo}{Theorem}
\newtheorem{coro}[theo]{Corollary}
\newtheorem{lemm}[theo]{Lemma}
\newtheorem{defi}[theo]{Definition}

\newcommand{\C}[1]{\mathcal{#1}}					
\newcommand{\setm}{\setminus}						
\newcommand{\eref}[1]{(\ref{#1})}					

\newcommand{\Prob}{\mathbb{P}}						
\newcommand{\given}{\,|\,}						
\newcommand{\normaldist}[2]{\C{N}\left(#1,{#2}^2\right)}		

\newcommand{\abs}[1]{\left\lvert #1 \right\rvert}			
\newcommand{\nel}[1]{\#(#1)}						
\newcommand{\norms}[2]{\,\left\Vert {#1} \right\Vert_{#2}\,}		
\newcommand{\NN}{\mathbb{N}}						

\newcommand{\graph}{\C{G}}						
\newcommand{\uedges}{\C{E}}						
\newcommand{\parents}[1]{\mathrm{par}(#1)}				
\newcommand{\children}[1]{\mathrm{ch}(#1)}				

\newcommand{\ue}[2]{\{#1,#2\}}						

\newcommand{\vars}{\C{V}}						
\newcommand{\CX}[1]{\C{X}_{#1}}						
\newcommand{\x}[1]{x_{#1}}						
\newcommand{\factors}{\C{F}}						
\newcommand{\FAC}{\psi}							
\newcommand{\fac}[1]{\FAC_{#1}}						

\newcommand{\nb}[1]{N_{#1}}						
\newcommand{\nbv}[1]{N_{#1}}						
\newcommand{\nbf}[1]{N_{#1}}						
\newcommand{\nbve}[2]{\nbv{#1}\setm #2}					
\newcommand{\nbfe}[2]{\nbf{#1}\setm #2}					
\newcommand{\del}[1]{{\partial{#1}}}					
\newcommand{\Del}[1]{{\Delta {#1}}}					

\newcommand{\xd}[1]{\x{\del{#1}}}

\newcommand{\bigFAC}{\Psi}						
\newcommand{\Fac}[1]{\bigFAC_{#1}}					

\newcommand{\Pm}[1]{\Prob^{\setm #1}}					
\newcommand{\Pmx}[1]{\Pm{#1}(\xd{#1})}					

\newcommand{\uFac}[1]{{\overline\Psi}_{#1}}				
\newcommand{\lFac}[1]{{\underline\Psi}_{#1}}				
\newcommand{\nor}{\C{N}}						
\newcommand{\parsum}{\C{Z}}						
\newcommand{\ep}{\mathrm{Ext}\,}					
\newcommand{\hull}{\mathrm{Hull}\,}					
\newcommand{\meas}[1]{\C{M}_{#1}}					
\newcommand{\fmeas}[1]{\C{Q}_{#1}}					
\newcommand{\nmeas}[1]{\C{P}_{#1}}					
\newcommand{\bbox}[3]{\C{B}_{#1}\left(#2,#3\right)}			
\newcommand{\sbbox}[1]{\C{B}\left(#1\right)}				
\newcommand{\snbbox}[1]{\C{B}\nor\left(#1\right)}			
\newcommand{\mebox}[2]{\C{B}_{#1 \to #2}}				
\newcommand{\bebox}[1]{\C{B}_{#1}}					

\ifthenelse{\isundefined{\techrep}}{\jmlrheading{vol}{year}{pages}{submitted}{published}{Joris M.\ Mooij and Hilbert J.\ Kappen}}{}

\ShortHeadings{Novel Bounds on Marginal Probabilities}{Mooij and Kappen}
\firstpageno{1}

\begin{document}

\title{Novel Bounds on Marginal Probabilities}

\author{\name Joris M.\ Mooij \email joris.mooij@tuebingen.mpg.de \\
\addr MPI for Biological Cybernetics, Dept. Sch\"olkopf \\
Spemannstra\ss e 38, 72076 T\"ubingen, Germany \\
\AND
\name Hilbert J. Kappen \email b.kappen@science.ru.nl \\
\addr Department of Biophysics \\
Radboud University Nijmegen \\
6525 EZ Nijmegen, The Netherlands}

\ifthenelse{\isundefined{\techrep}}{\editor{Editor}}{\editor{}}

\maketitle

\begin{abstract}%
We derive two related novel bounds on single-variable marginal probability 
distributions in factor graphs with discrete variables. The first 
method propagates bounds over a subtree of the factor graph rooted in
the variable, and the second method propagates bounds over the
self-avoiding walk tree starting at the variable. By 
construction, both methods not only bound the exact marginal probability 
distribution of a variable, but also its approximate Belief Propagation 
marginal (``belief''). Thus, apart from providing a practical means to 
calculate bounds on marginals, our contribution also lies in an increased
understanding of the error made by Belief Propagation. Empirically, we show 
that our bounds often outperform existing bounds in terms of accuracy 
and/or computation time. We also show that our bounds can yield nontrivial 
results for medical diagnosis inference problems.
\end{abstract}


\section{Introduction}

Graphical models are used in many different fields. A fundamental problem in
the application of graphical models is that exact inference is NP-hard
\citep{Cooper90}. In recent years, much research has focused on approximate
inference techniques, such as sampling methods and deterministic approximation
methods, \eg Belief Propagation (BP) \citep{Pearl88}.  Although the
approximations obtained by these methods can be very accurate, there are only
few guarantees on the error of the approximation, and often it is not known
(without comparing with the exact solution) how accurate an approximate result
is. Thus it is desirable to calculate, in addition to the approximate results,
tight bounds on the approximation error. Existing methods to calculate bounds 
on marginals include \citep{Tatikonda03,LeisinkKappen03,TagaMase06b,Ihler07}.
Also, upper bounds on the partition sum, \eg \citep{JaakkolaJordan96,WainwrightJaakkolaWillsky05},
can be combined with lower bounds on the partition sum, such as the well-known
mean field bound or higher-order lower bounds \citep{LeisinkKappen01}, to 
obtain bounds on marginals.

In this article, we derive novel bounds on exact single-variable marginals in
factor graphs. The original motivation for this work was to better
understand and quantify the BP error. This has lead to bounds which are at the
same time bounds for the exact single-variable marginals as well as for the BP
beliefs. A particularly nice feature of the bounds is that their computational
cost is relatively low, provided that the number of possible values of each
variable in the
factor graph is small. Unfortunately, the computation time is exponential
in the number of possible values of the variables, which limits application to
factor graphs in which each variable has a low number of possible values.
On these factor graphs however, our bounds perform exceedingly well and we
show empirically that they outperform the state-of-the-art in a variety of
factor graphs, including real-world problems arising in medical diagnosis.

This article is organized as follows.
In the next section, we derive our novel bounds. In Section \ref{sec5:related_work},
we discuss related work. In Section \ref{sec5:experiments} we present experimental results. 
We conclude with conclusions and a discussion in Section \ref{sec5:discussion}.

\section{Theory}\label{sec5:theory}

In this work, we consider graphical models such as Markov random fields and
Bayesian networks. We use the unifying factor graph representation
\citep{KschischangFreyLoeliger01}. In the first subsection, we introduce our
notation and some basic definitions concerning factor graphs. Then, we shortly
remind the reader of some basic facts about convexity. After that, we introduce
some notation and concepts for measures on subsets of variables. We proceed
with a subsection that considers the interplay between convexity and the
operations of normalization and multiplication.  In the next subsection, we
introduce ``(smallest bounding) boxes'' that will be used to describe sets of
measures in a convenient way. Then, we formulate the basic lemma that will be
used to obtain bounds on marginals. We illustrate the basic lemma with two
simple examples. Then we formulate our first result, an algorithm for
propagating boxes over a subtree of the factor graph, which results in
a bound on the marginal of the root variable of the subtree. In the last
subsection, we show how one can go deeper into the computation tree and derive
our second result, an algorithm for propagating boxes over self-avoiding walk
trees. The result of that algorithm is a bound on the marginal of the root
variable (starting point) of the self-avoiding walk tree. For the special case
where all factors in the factor graph depend on two variables at most 
(``pairwise interactions''), our
first result is equivalent to a truncation of the second one. This is not true
for higher-order interactions, however.

\subsection{Factor graphs}

Let $\vars := \{1,\dots,N\}$ and consider $N$ discrete random variables 
$(\x{i})_{i\in\vars}$. Each variable $\x{i}$ takes values in a discrete 
domain $\CX{i}$. We will frequently use the following multi-index notation. Let
$A = \{i_1, i_2, \dots, i_m\} \subseteq \vars$ with $i_1 < i_2 < \dots i_m$.
We write $\CX{A} := \CX{i_1} \times \CX{i_2} \times \dots \times \CX{i_m}$ 
and for any family $(Y_i)_{i\in B}$ with $A \subseteq B \subseteq \vars$,
we write $Y_A := (Y_{i_1}, Y_{i_2}, \dots, Y_{i_m})$. 

We consider a probability distribution over 
$x = (\x{1},\dots,\x{N}) \in \CX{\vars}$ that can be written as a product of factors 
(also called ``interactions'') $(\fac{I})_{I\in\factors}$:
  \begin{equation}\label{eq:prob_dist}
  \Prob(x) = \frac{1}{Z} \prod_{I\in\factors} \fac{I}(\x{\nbf{I}}), \qquad Z = \sum_{x\in\CX{\vars}} \prod_{I \in \factors} \fac{I}(\x{\nbf{I}}).
  \end{equation}
For each factor index $I \in \factors$, there is an associated subset 
$\nbf{I} \subseteq\vars$ of variable indices and the factor $\fac{I}$ 
is a nonnegative function $\fac{I} : \CX{\nbf{I}} \to [0,\infty)$.
For a Bayesian network, the factors are (conditional) probability
tables. In case of Markov random fields, the factors are often called
potentials.\footnote{Not to be confused with statistical physics 
terminology, where ``potential'' refers to $-\frac{1}{\beta}\log \fac{I}$ 
instead, with $\beta$ the inverse temperature.} 
In the following, we will use lowercase for variable indices and uppercase for
factor indices. 

In general, the normalizing constant $Z$ is not known and exact computation of
$Z$ is infeasible, due to the fact that the number of terms to be summed is
exponential in $N$. Similarly, computing marginal distributions $\Prob(\x{A})$ 
for subsets of variables $A \subseteq \vars$ is intractable in general.
In this article, we focus on the task of obtaining rigorous bounds on single-variable
marginals $\Prob(\x{i}) = \sum_{\x{\vars\setm\{i\}}}\Prob(x)$.

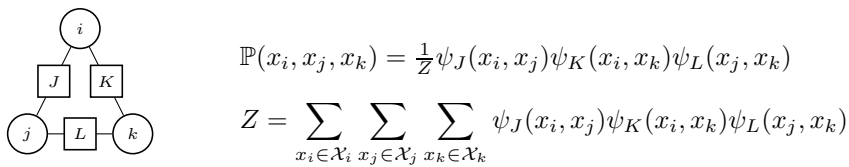
\begin{figure}
  \centering
  \begin{tikzpicture}
    \tikzstyle{var}=[circle,draw=black,fill=white,semithick,minimum size=15pt]
    \tikzstyle{fac}=[rectangle,draw=black,fill=white,semithick,minimum size=12pt]
    \small
    \begin{scope}[scale=0.7]
      \tiny
      \node (i) at (0,2)  [var] {$i$};
      \node (j) at (-1,0) [var] {$j$};
      \node (k) at (1,0)  [var] {$k$};
      \path (i) edge node[pos=0.5] (J) [fac] {$J$} (j);
      \path (i) edge node[pos=0.5] (K) [fac] {$K$} (k); 
      \path (j) edge node[pos=0.5] (L) [fac] {$L$} (k);
    \end{scope}
    \small
    \node at (2,1) [anchor=west] {$\Prob(\x{i},\x{j},\x{k}) = \frac{1}{Z} \fac{J}(\x{i},\x{j}) \fac{K}(\x{i},\x{k}) \fac{L}(\x{j},\x{k})$};
    \node at (2,0) [anchor=west] {$Z = \displaystyle\sum_{\x{i}\in\CX{i}} \sum_{\x{j}\in\CX{j}} \sum_{\x{k}\in\CX{k}} \fac{J}(\x{i},\x{j}) \fac{K}(\x{i},\x{k}) \fac{L}(\x{j},\x{k})$};
  \end{tikzpicture} 
\caption{\label{fig:factorgraph_simple}Example of a factor graph with three variable nodes ($i,j,k$), represented
as circles, and three factor nodes ($J,K,L$), represented as rectangles. The corresponding random variables are $\x{i},\x{j},\x{k}$;
the corresponding factors are $\fac{J} : \CX{i} \times \CX{j} \to [0,\infty)$,
$\fac{K} : \CX{i} \times \CX{k} \to [0,\infty)$ and $\fac{L} : \CX{j} \times \CX{k} \to [0,\infty)$.
The corresponding probability distribution $\Prob(\x{})$ is written out on the right.}
\end{figure}

We can represent the structure of the probability distribution
\eref{eq:prob_dist} using a \emph{factor graph} $(\vars,\factors,\uedges)$. 
This is a bipartite graph, consisting of \emph{variable nodes} $i \in \vars$, 
\emph{factor nodes} $I \in \factors$, and \emph{edges} $e \in \uedges$, 
with an edge $\ue{i}{I}$ between $i\in\vars$ and $I\in\factors$ if and only if the factor $\fac{I}$ 
depends on $\x{i}$ (\ie if $i \in \nbf{I}$). We will
represent factor nodes visually as rectangles and variable nodes as circles.
Figure \ref{fig:factorgraph_simple} shows a simple example of a factor graph
and the corresponding probability distribution. The set of neighbors of 
a factor node $I$ is precisely $\nbf{I}$; similarly, we denote the set of
neighbors of a variable node $i$ by $\nbv{i} := \{I \in \factors: i \in \nbf{I}\}$.
Further, we define for each variable 
$i \in \vars$ the set $\Del{i} := \bigcup \nbv{i}$ consisting of all variables 
that appear in some factor in which variable $i$ participates, and the set 
$\del{i} := \Del{i} \setm \{i\}$, the \emph{Markov blanket} of $i$. 

We will assume throughout this article that the factor
graph corresponding to \eref{eq:prob_dist} is connected.
Furthermore, we will assume that
\begin{equation*}
\forall{I \in \factors} \ \forall{i\in \nbf{I}} \ \forall {\x{\nbfe{I}{\{i\}}} \in \CX{\nbfe{I}{\{i\}}}}: \sum_{\x{i} \in \CX{i}} \fac{I}(\x{i},\x{\nbfe{I}{\{i\}}}) > 0.
\end{equation*}
This will prevent technical problems regarding normalization later on.\footnote{This
condition ensures that if one runs Belief Propagation on the factor graph,
the messages will always remain nonzero, provided that the initial messages
are nonzero.}

One final remark concerning notation: we will sometimes abbreviate $\{i\}$
as $i$ if no confusion can arise.

\subsection{Convexity}

Let $V$ be a real vector space. For $T$ elements $(v_t)_{t=1,\dots,T}$ of $V$
and $T$ nonnegative numbers $(\lambda_t)_{t=1,\dots,T}$ with $\sum_{t=1}^T
\lambda_t = 1$, we call $\sum_{t=1}^T \lambda_t v_t$ a \emph{convex
combination} of the $(v_t)_{t=1,\dots,T}$ with weights $(\lambda_t)_{t=1,\dots,T}$.  
A subset $X \subseteq V$ is called
\emph{convex} if for all $x_1,x_2 \in X$ and all $\lambda \in [0,1]$, the
convex combination $\lambda x_1 + (1-\lambda)x_2 \in X$.  An
\emph{extreme point} of a convex set $X$ is an element $x \in X$ which cannot
be written as a (nontrivial) convex combination of two different points in $X$.
In other words, $x \in X$ is an extreme point of $X$ if and only if for all
$\lambda \in (0,1)$ and all $x_1,x_2 \in X$, $x = \lambda x_1 + (1-\lambda)
x_2$ implies $x_1 = x_2$.  We denote the set of extreme points of a convex set
$X$ by $\ep(X)$.  For a subset $Y$ of the vector space $V$, we define the
\emph{convex hull of $Y$} to be the smallest convex set $X \subseteq V$ with
$Y \subseteq X$; we denote the convex hull of $Y$ as $\hull(Y)$. 

\subsection{Measures and operators}

For $A \subseteq \vars$, define $\meas{A} := [0,\infty)^{\CX{A}}$, 
\ie $\meas{A}$ is the set of nonnegative functions on $\CX{A}$. $\meas{A}$ can 
be identified with the set of finite measures on $\CX{A}$. We will simply call the 
elements of $\meas{A}$ ``measures on $A$''.
We also define $\meas{A}^* := \meas{A} \setm \{0\}$.
We will denote $\meas{} := \bigcup_{A \subseteq \vars} \meas{A}$ and
$\meas{}^* := \bigcup_{A \subseteq \vars} \meas{A}^*$.

Adding two measures $\Psi,\Phi\in\meas{A}$ results in the measure $\Psi+\Phi$
in $\meas{A}$.
For $A, B \subseteq \vars$, we can multiply an element of $\meas{A}$ with
an element of $\meas{B}$ to obtain an element of $\meas{A\cup B}$; a special
case is multiplication with a scalar. 
Note that there is a natural embedding of $\meas{A}$ in $\meas{B}$ for $A \subseteq B \subseteq \vars$
obtained by multiplying an element $\Psi\in\meas{A}$ by $\mathbf{1}_{B\setm A}\in\meas{B\setm A}$,
the constant function with value 1 on $\CX{B\setm A}$.
Another important operation is the partial
summation: given $A \subseteq B \subseteq \vars$ and $\Psi\in\meas{B}$, define
$\sum_{x_A} \Psi$ to be the measure in $\meas{B\setm A}$ that satisfies
 $$\left(\sum_{x_A} \Psi\right)(x_{B\setm A}) = \sum_{x_A\in\CX{A}} \Psi(x_A, x_{B\setm A}) \qquad \forall x_{B\setm A} \in \CX{B\setm A}.$$
Also, defining $A'=B\setm A$, we will sometimes write this measure as $\sum_{\x{\setm A'}} \Psi$,
which is an abbreviation of $\sum_{\x{B\setm A'}} \Psi$. This notation does not make
explicit which variables are summed over (which depends on the measure that is being partially
summed), although it shows which variables remain after summation.

In the following, we will implicitly define operations on \emph{sets} of
measures by applying the operation on elements of these sets and taking
the set of the resulting measures; e.g., if we have two subsets $\Xi_A\subseteq\meas{A}$
and $\Xi_B\subseteq\meas{B}$ for $A,B\subseteq \vars$, we define the product
of the sets $\Xi_A$ and $\Xi_B$ to be the set of the products of elements of
$\Xi_A$ and $\Xi_B$, \ie $\Xi_A \Xi_B := \{\Psi_A \Psi_B : \Psi_A \in \Xi_A, \Psi_B \in \Xi_B\}$.

\begin{figure}
\centering
\begin{tikzpicture}
  \fill[fill=gray!10] (0,0) rectangle +(4,3);
  \draw [->,thick] (0,0) -- (4,0);
  \draw [->,thick] (0,0) -- (0,3);
  \node at (4,0) [anchor=north] {\tiny$\x{i}=+1$};
  \node at (0,3) [anchor=south] {\tiny$\x{i}=-1$};
  \node at (4,3) [anchor=north east] {$\meas{i}$};
  \draw [draw=black!80,thick] (2,0) -- (0,2);
  \node at (0,2) [anchor=east] {\tiny1};
  \node at (2,0) [anchor=north] {\tiny1};
  \node at (0.5,1.5) [anchor=south,xshift=0.1cm] {$\nmeas{i}$};
  \draw [draw=gray,thin] (3,2) -- (1.2,0.8) -- (0,0);
  \node (Psi) at (3,2) [anchor=west] {$\Fac{}$};
  \node (NPsi) at (1.2,0.8) [right=5pt,anchor=west] {$\nor \Fac{}$};
  \filldraw (1.2,0.8) circle (1.5pt);
  \filldraw (3,2) circle (1.5pt);
\end{tikzpicture}
\caption{\label{fig:measures_binary}Illustration of some concepts in the simple case of a 
binary random variable $\x{i} \in \CX{i} = \{\pm 1\}$ and the subset $A = \{i\}$.
A measure $\Fac{} \in \meas{i}$ can be identified with a point in the quarter
plane as indicated in the Figure. A normalized measure can be
obtained by scaling $\Fac{}$; the result $\nor\Fac{}$ is contained in the simplex $\nmeas{i}$,
a lower-dimensional submanifold of $\meas{i}$.}
\end{figure}
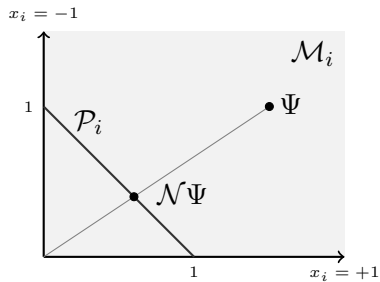

In Figure \ref{fig:measures_binary}, the simple case of a binary random 
variable $\x{i}$ and the subset $A = \{i\}$ is illustrated. Note that
in this case, a measure $\Fac{} \in \meas{i}$ can be identified with a
point in the quarter plane $[0,\infty) \times [0,\infty)$.

We will define $\fmeas{A}$ to be the set of completely factorized measures on $A$, \ie
$$\fmeas{A} := \prod_{a \in A} \meas{\{a\}} = \left\{\prod_{a\in A} \Fac{a} : \Fac{a} \in \meas{\{a\}} \text{ for each } a \in A \right\}.$$
Note that $\meas{A}$ is the convex hull of $\fmeas{A}$. Indeed, we can write each measure
$\Psi\in\meas{A}$ as a convex combination of measures in $\fmeas{A}$; let 
$Z := \sum_{\x{A}} \Psi$ and note that
  $$\Psi(x) = \sum_{y\in \CX{A}} \frac{\Psi(y)}{Z} \left(Z \delta_y(x)\right)\qquad \forall x \in \CX{A}.$$
For any $y \in \CX{A}$, the Kronecker delta function $\delta_y \in \meas{A}$ (which is 1 if its
argument is equal to $y$ and 0 otherwise) is
an element of $\fmeas{A}$ because $\delta_y(x) = \prod_{a\in A} \delta_{y_a}(x_a)$.
We denote $\fmeas{A}^* := \fmeas{A} \setm \{0\}$.

We define the \emph{partition sum operator} $\parsum : \meas{} \to [0,\infty)$ which calculates
the partition sum (normalization constant) of a measure, \ie
\begin{equation*}
\parsum \Fac{} := \sum_{\x{A}\in\CX{A}} \Fac{}(\x{A})\qquad \text{for $\Fac{} \in \meas{A}$, $A\subseteq\vars$}.
\end{equation*}
We denote with $\nmeas{A}$ the set of probability measures on $A$, \ie 
$\nmeas{A} = \{\Fac{} \in \meas{A} : \parsum \Fac{} = 1\}$, and define
$\nmeas{} := \bigcup_{A \subseteq \vars} \nmeas{A}$.
The set $\nmeas{A}$ is called a \emph{simplex} (see also Figure \ref{fig:measures_binary}).
Note that a simplex is convex; the simplex $\nmeas{A}$ has precisely $\nel{\CX{A}}$ extreme
points, each of which corresponds to putting all probability mass on one of the possible values
of $\x{A}$.

Define the \emph{normalization operator} $\nor : \meas{}^* \to \nmeas{}$ 
which normalizes a measure, \ie
\begin{equation*}
\nor \Fac{} := \frac{1}{\parsum \Fac{}} \Fac{} \qquad \text{for $\Fac{} \in \meas{}^*$}.
\end{equation*}
Note that $\parsum \circ \nor = 1$. Figure \ref{fig:measures_binary} illustrates the normalization
of a measure in a simple case.

\subsection{Convex sets of measures}

To calculate marginals of subsets of variables in some factor graph, several
operations performed on measures are relevant: normalization, taking products
of measures, and summing over subsets of variables. In this section we study
the interplay between convexity and these operations; this will turn out to be
useful later on, because our bounds make use of convex sets of measures that
are propagated over the factor graph.

\begin{figure}
\centering
\begin{tikzpicture}
  \fill[fill=gray!10] (0,0) rectangle +(4,3.5);
  \draw [->,thick] (0,0) -- (4,0);
  \draw [->,thick] (0,0) -- (0,3.5);
  \node at (4,0) [anchor=north] {\tiny$\x{i}=+1$};
  \node at (0,3.5) [anchor=south] {\tiny$\x{i}=-1$};
  \node at (4,3.5) [anchor=north east] {$\meas{i}$};
  \draw [draw=black!80,semithick] (2,0) -- (0,2);
  \node at (0,2) [anchor=east] {\tiny1};
  \node at (2,0) [anchor=north] {\tiny1};
  \node at (0.2,1.8) [anchor=south,xshift=0.1cm] {$\nmeas{i}$};
  \draw [very thick] (0.6,1.4) -- (1.5,0.5);
  \filldraw [fill=gray,draw=gray] (1.2,2.8) -- (1.5,1.5) -- (3.0,1.0);
  \draw [draw=gray,thin] (1.2,2.8) -- (0,0);
  \draw [draw=gray,thin] (1.5,1.5) -- (0,0);
  \draw [draw=gray,thin] (3.0,1.0) -- (0,0);
  \node (NPsi) at (1.2,2.8) [right=2pt,anchor=east] {\tiny$\xi_1$};
  \node (NPsi) at (1.5,1.5) [anchor=east,yshift=1pt] {\tiny$\xi_2$};
  \node (NPsi) at (3.0,1.0) [anchor=west] {\tiny$\xi_3$};
  \node (NPsi) at (0.6,1.4) [right=2pt,anchor=east] {\tiny$\nor \xi_1$};
  \node (NPsi) at (1.0,1.0) [anchor=west] {\tiny$\nor \xi_2$};
  \node (NPsi) at (1.5,0.5) [anchor=east,yshift=2pt] {\tiny$\nor \xi_3$};
  \filldraw (0.6,1.4) circle (1.5pt);
  \filldraw (1.0,1.0) circle (1.5pt);
  \filldraw (1.5,0.5) circle (1.5pt);
  \filldraw (1.2,2.8) circle (1.5pt);
  \filldraw (1.5,1.5) circle (1.5pt);
  \filldraw (3.0,1.0) circle (1.5pt);
\end{tikzpicture}
\caption{\label{fig:normalizing_convexity}Any convex combination of $\nor\xi_1$, $\nor\xi_2$ and $\nor\xi_3$
can be written as a normalized convex combination of $\xi_1$, $\xi_2$ and $\xi_3$. Vice versa, normalizing a
convex combination of $\xi_1,\xi_2$ and $\xi_3$ yields a convex combination of $\nor\xi_1, \nor\xi_2$ and $\nor\xi_3$.}
\end{figure}
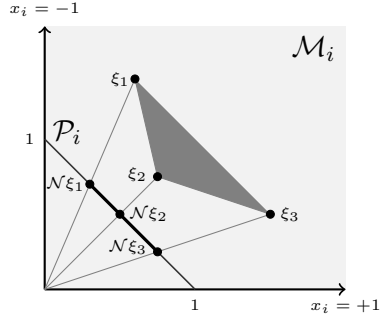

The interplay between normalization and convexity is described by the following 
Lemma, which is illustrated in Figure \ref{fig:normalizing_convexity}.

\begin{lemm}\label{lemm:normalizing_preserves_convexity}
Let $A \subseteq \vars$, $T\in\NN^*$ and let
$(\xi_t)_{t=1,\dots,T}$ be elements of $\meas{A}^*$.
Each convex combination of the normalized measures
$(\nor \xi_t)_{t=1,\dots,T}$ can be written as a normalized
convex combination of the measures $(\xi_t)_{t=1,\dots,T}$
(which has different weights in general), and vice versa.
\end{lemm}
\begin{proof}
Let $(\lambda_t)_{t=1,\dots,T}$ be nonnegative numbers
with $\sum_{t=1}^T \lambda_t = 1$. Then
\begin{equation*}
\parsum \left(\sum_{t=1}^T \lambda_t \nor \xi_t\right) = 
\sum_{t=1}^T \lambda_t \parsum \nor \xi_t = 1,
\end{equation*}
therefore
\begin{equation*}
\sum_{t=1}^T \lambda_t (\nor \xi_t)
 = \nor \left( \sum_{t=1}^T \lambda_t (\nor \xi_t) \right) 
 = \nor \left( \sum_{t=1}^T \frac{\lambda_t}{\parsum \xi_t} \xi_t \right) 
 = \nor \left( \sum_{t=1}^T \frac{\frac{\lambda_t}{\parsum \xi_t}}{\sum_{s=1}^T \frac{\lambda_s}{\parsum \xi_s}} \ \xi_t \right),
\end{equation*}
which is the result of applying the normalization operator
to a convex combination of the elements $(\xi_t)_{t=1,\dots,T}$.

Vice versa, let $(\mu_t)_{t=1,\dots,T}$ be nonnegative numbers with $\sum_{t=1}^T \mu_t = 1$. Then
$$\nor \left( \sum_{t=1}^T \mu_t \xi_t \right) = \sum_{t=1}^T \frac{\mu_t}{Z} \xi_t$$
where
$$Z := \parsum \left( \sum_{t=1}^T \mu_t \xi_t \right) = \sum_{t=1}^T \mu_t \parsum \xi_t = \sum_{t=1}^T \mu_t Z_t$$
where we defined $Z_t := \parsum \xi_t$ for all $t=1,\dots,T$. Thus
$$\nor \left( \sum_{t=1}^T \mu_t \xi_t \right) = \sum_{t=1}^T \frac{\mu_t}{\sum_{s=1}^T \mu_s Z_s} \xi_t
= \sum_{t=1}^T \frac{\mu_t Z_t}{\sum_{s=1}^T \mu_s Z_s} \nor \xi_t,$$
which is a convex combination of the normalized measures $(\nor \xi_t)_{t=1,\dots,T}$.
\end{proof}

The following lemma concerns the interplay between convexity
and taking products; it says that if we take the product of convex
sets of measures on different spaces, the resulting set is
contained in the convex hull of the product of the extreme
points of the convex sets. We have not made a picture corresponding
to this lemma because the simplest nontrivial case would require at least 
four dimensions.
\begin{lemm}\label{lemm:product_convex_sets}
Let $T \in \NN^*$ and $(A_t)_{t=1,\dots,T}$
be a family of mutually disjoint subsets of $\vars$. For each $t=1,\dots,T$,
let $\Xi_t \subseteq \meas{A_t}$ be convex with a finite number of extreme points. Then:
\begin{equation*}
\prod_{t=1}^T \Xi_t \subseteq \hull \left(\prod_{t=1}^T \ep \Xi_t \right),
\end{equation*}
\end{lemm}
\begin{proof}
Let $\Fac{t} \in \Xi_t$ for each $t=1,\dots,T$.
For each $t$, $\Fac{t}$ can be written as a convex combination 
\begin{equation*}
\Fac{t} = \sum_{\xi_t \in \ep(\Xi_t)} \lambda_{t;\xi_t} \xi_t,
\qquad \sum_{\xi_t \in \ep(\Xi_t)} \lambda_{t;\xi_t} = 1,
\qquad \forall\ {\xi_t \in \ep(\Xi_t)}: \lambda_{t;\xi_t} \ge 0.
\end{equation*}
Therefore the product $\prod_{t=1}^T \Fac{t}$ is also a convex combination:
\begin{align*}
\prod_{t=1}^T \Fac{t} 
& = \prod_{t=1}^T \left( \sum_{\xi_t \in \ep(\Xi_t)} \lambda_{t;\xi_t} \xi_t \right) \\
& = \sum_{\xi_1 \in \ep(\Xi_1)} \sum_{\xi_2\in\ep(\Xi_2)} \dots \sum_{\xi_T\in\ep(\Xi_T)} \left( \prod_{t=1}^T \lambda_{t;\xi_t} \right) \left( \prod_{t=1}^T \xi_t \right) \\
& \in \hull \left( \prod_{t=1}^T \ep \Xi_t \right).
\end{align*}
\end{proof}

\subsection{Boxes and smallest bounding boxes}

In this subsection, we define ``(smallest bounding) boxes'', certain convex sets of
measures that will play a central role in our bounds, and study some of
their properties.

\begin{defi}
Let $A \subseteq \vars$.
For $\lFac{} \in \meas{A}$ and $\uFac{} \in \meas{A}$
with $\lFac{} \le \uFac{}$, we define the \emph{box} between
the lower bound $\lFac{}$ and the upper bound $\uFac{}$ by
\begin{equation*}
\bbox{A}{\lFac{}}{\uFac{}} := \{\Fac{} \in \meas{A} : \lFac{} \le \Fac{} \le \uFac{}\}.
\end{equation*}
\end{defi}
The inequalities should be interpreted pointwise, \eg $\lFac{} \le \Fac{}$ means
$\lFac{}(x) \le \Fac{}(x)$ for all $x \in \CX{A}$.
Note that a box is convex; indeed, its extreme points are the
``corners'' of which there are $2^{\nel{\CX{A}}}$.

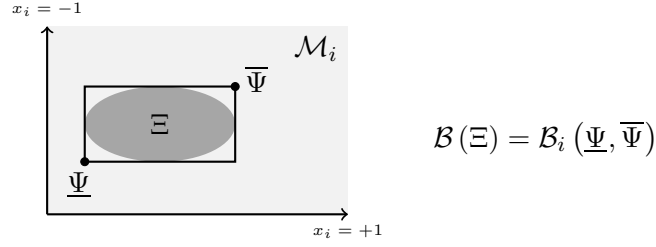
\begin{figure} 
\centering
\begin{tikzpicture} 
  \fill[fill=gray!10] (0,0) rectangle +(4,2.5);
  \draw [->,thick] (0,0) -- (4,0);
  \draw [->,thick] (0,0) -- (0,2.5);
  \node at (4,0)   [anchor=north] {\tiny$\x{i}=+1$};
  \node at (0,2.5) [anchor=south] {\tiny$\x{i}=-1$};
  \node at (4,2.5) [anchor=north east] {$\meas{i}$};
  \fill[fill=gray!70] (1.5,1.2) ellipse (1cm and 0.5cm);
  \node (Xi) at (1.5,1.2) {$\Xi$};
  \draw[draw=black,thick] (0.5,0.7) rectangle +(2,1);
  \filldraw[fill=black,draw=black] (0.5,0.7) circle (1.5pt);
  \filldraw[fill=black,draw=black] (2.5,1.7) circle (1.5pt);
  \node (PsiLo) at (0.5,0.7) [anchor=north,xshift=-2pt] {$\lFac{}$};
  \node (PsiHi) at (2.5,1.7) [anchor=west,yshift=2pt] {$\uFac{}$};
  \node (formula) at (5,1) [anchor=west] {$\sbbox{\Xi} = \bbox{i}{\lFac{}}{\uFac{}}$};
\end{tikzpicture}
\caption{\label{fig:bounding_boxes}The smallest bounding box $\sbbox{\Xi}$ for $\Xi$ is given by
the box $\bbox{i}{\lFac{}}{\uFac{}}$ with lower bound $\lFac{}$ and upper bound $\uFac{}$.}
\end{figure}

\begin{defi}
Let $A \subseteq \vars$ and $\Xi \subseteq \meas{A}$ be bounded (\ie $\Xi \le \Fac{}$ for some $\Fac{} \in \meas{A}$).
The \emph{smallest bounding box} for $\Xi$ is defined as
$\sbbox{\Xi} := \bbox{A}{\lFac{}}{\uFac{}}$,
where $\lFac{}, \uFac{} \in \meas{A}$ are given by
\begin{equation*}
\begin{array}{lll}
\lFac{}(\x{A}) &\! := \inf \{ \Fac{}(\x{A}) : \Fac{} \in \Xi \} & \qquad \forall \x{a}\in\CX{A},\\
\uFac{}(\x{A}) &\! := \sup \{ \Fac{}(\x{A}) : \Fac{} \in \Xi \} & \qquad \forall \x{a}\in\CX{A}.
\end{array}
\end{equation*}
\end{defi}
Figure \ref{fig:bounding_boxes} illustrates this concept.
Note that $\sbbox{\Xi} = \sbbox{\hull(\Xi)}$. Therefore, if
$\Xi$ is convex, the smallest bounding box for $\Xi$ depends only
on the extreme points $\ep(\Xi)$, \ie  $\sbbox{\Xi} = \sbbox{\ep(\Xi)}$.

\begin{figure}
\centering
\begin{tikzpicture}
  \fill[fill=gray!10] (0,0) rectangle +(4,4);
  \draw [->,thick] (0,0) -- (4,0);
  \draw [->,thick] (0,0) -- (0,4);
  \node at (4,0) [anchor=north] {\tiny$\x{i}=+1$};
  \node at (0,4) [anchor=south] {\tiny$\x{i}=-1$};
  \node at (4,4) [anchor=north east] {$\meas{i}$};
  \draw (1.6,0.8) rectangle (3.3,2.2);          
  \draw (0.8,1.8) rectangle (2.7,2.6);          
  \filldraw[fill=black,draw=black] (1.6,0.8) circle (1.5pt);
  \filldraw[fill=black,draw=black] (3.3,2.2) circle (1.5pt);
  \filldraw[fill=black,draw=black] (0.8,1.8) circle (1.5pt);
  \filldraw[fill=black,draw=black] (2.7,2.6) circle (1.5pt);
  \node (PsiLo1) at (1.6,0.8) [anchor=east] {$\lFac{1}$};
  \node (PsiHi1) at (3.3,2.2) [anchor=west,yshift=-2pt] {$\uFac{1}$};
  \node (PsiLo2) at (0.8,1.8) [anchor=east] {$\lFac{2}$};
  \node (PsiHi2) at (2.7,2.6) [anchor=west,yshift=2pt] {$\uFac{2}$};
\begin{scope}[xshift=5.5cm]
  \fill[fill=gray!10] (0,0) rectangle +(5,4);
  \draw [->,thick] (0,0) -- (5,0);
  \draw [->,thick] (0,0) -- (0,4);
  \node at (5,0) [anchor=north] {\tiny$\x{i}=+1$};
  \node at (0,4) [anchor=south] {\tiny$\x{i}=-1$};
  \node at (5,4) [anchor=north east] {$\meas{i}$};
  \draw (0.64,0.72) rectangle (4.46,2.86);
  \filldraw[fill=black,draw=black] (0.64,0.72) circle (1.5pt);
  \filldraw[fill=black,draw=black] (4.46,2.86) circle (1.5pt);
  \node (PsiLo12) at (0.64,0.72) [anchor=north] {$\lFac{1}\lFac{2}$};
  \node (PsiHi12) at (4.46,2.86) [anchor=south] {$\uFac{1}\uFac{2}$};
  \node (bbox) at (2.55,1.79) {\small $\bbox{i}{\lFac{1}}{\uFac{1}} \bbox{i}{\lFac{2}}{\uFac{2}}$};
\end{scope}
\end{tikzpicture}
\caption{\label{fig:product_of_bbs}Multiplication of two boxes on the same variable set $A = \{i\}$.}
\end{figure}
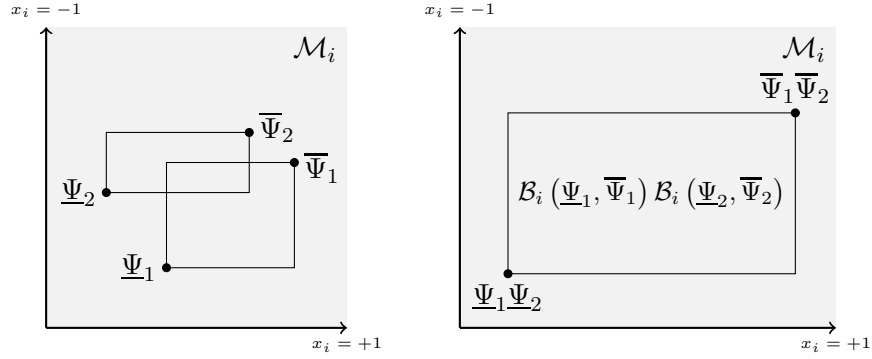

The product of several boxes on the same subset $A$ of variables can
be easily calculated as follows (see also Figure \ref{fig:product_of_bbs}).
\begin{lemm}\label{lemm:product_boxes_same_varset}
Let $A \subseteq \vars$, $T \in \NN^*$ and for each $t = 1, \dots, T$,
let $\lFac{t}, \uFac{t} \in \meas{A}$ such that $\lFac{t} \le \uFac{t}$.
Then
\begin{equation*}
\prod_{t=1}^T \bbox{A}{\lFac{t}}{\uFac{t}} = \bbox{A}{\prod_{t=1}^T \lFac{t}}{\prod_{t=1}^T \uFac{t}},
\end{equation*}
\ie the product of the boxes is again a box, with as lower bound the
product of the lower bounds of the boxes and as upper bound the
product of the upper bounds of the boxes.
\end{lemm}
\begin{proof}
We prove the case $T=2$; the general case follows by induction.
We show that
$$\bbox{A}{\lFac{1}}{\uFac{1}} \bbox{A}{\lFac{2}}{\uFac{2}} = \bbox{A}{\lFac{1}\,\lFac{2}}{\uFac{1}\,\uFac{2}}.$$
That is, for $\Phi \in \meas{A}$ we have to show that
\begin{equation*}
\lFac{1}(x)\lFac{2}(x) \le \Phi(x) \le \uFac{1}(x)\uFac{2}(x) \qquad \forall x \in \CX{A}
\end{equation*}
if and only if there exist $\Phi_1, \Phi_2 \in \meas{A}$ such that:
\begin{align*}
& \Phi(x) = \Phi_1(x) \Phi_2(x) & \qquad \forall x \in \CX{A}; \\
& \lFac{1}(x) \le \Phi_1(x) \le \uFac{1}(x) & \qquad \forall x \in \CX{A};\\
& \lFac{2}(x) \le \Phi_2(x) \le \uFac{2}(x) & \qquad \forall x \in \CX{A}.
\end{align*}
Note that the problem ``decouples'' for the various possible values of $x \in \CX{A}$
so that we can treat each component (indexed by $x \in \CX{A}$) seperately. That is, the problem
reduces to showing that
\begin{equation*}
[a, b] \cdot [c, d] = [ac, bd]
\end{equation*}
for $0 \le a \le b$ and $0 \le c \le d$ (take $a = \lFac{1}(x)$, $b = \uFac{1}(x)$,
$c = \lFac{2}(x)$ and $d = \uFac{2}(x)$). In other words, we have to show that $y \in [ac,bd]$
if and only if there exist $y_1 \in [a,b]$, $y_2 \in [c,d]$ with $y = y_1 y_2$.
For the less trivial part of this assertion, it is easily verified that choosing $y_1$
and $y_2$ according to the following table:
\begin{center}
\begin{tabular}{l|ll}
Condition           & $y_1$         & $y_2$ \\
\hline\\[-10pt]
$bc \le y, b > 0$   & $b$           & $\frac{y}{b}$ \\
$b = 0$             & $0$           & $c$ \\
$bc \ge y, c > 0$   & $\frac{y}{c}$ & $c$ \\
$bc \ge y, c = 0$   & $b$           & $0$
\end{tabular}
\end{center}
does the job.
\end{proof}

In general, the product of several boxes is not a box itself. Indeed, let $i,j \in \vars$
be two different variable indices. Then $\bbox{i}{\lFac{i}}{\uFac{i}} \bbox{j}{\lFac{j}}{\uFac{j}}$
contains only factorizing measures, whereas $\bbox{\{i,j\}}{\lFac{i}\lFac{j}}{\uFac{i}\uFac{j}}$
is not a subset of $\fmeas{\{i,j\}}$ in general. 
However, we do have the following identity:
\begin{lemm}
Let $T \in \NN^*$ and for each $t = 1, \dots, T$,
let $A_t \subseteq \vars$ and $\lFac{t}, \uFac{t} \in \meas{A_t}$ 
such that $\lFac{t} \le \uFac{t}$. Then
\begin{equation*}
\sbbox{\prod_{t=1}^T \bbox{A_t}{\lFac{t}}{\uFac{t}}} = \bbox{(\bigcup_{t=1}^T A_t)}{\prod_{t=1}^T \lFac{t}}{\prod_{t=1}^T \uFac{t}}.
\end{equation*}
\end{lemm}
\begin{proof}
Straightforward, using the definitions.
\end{proof}

\subsection{The basic lemma}

\begin{figure}
\centering
\begin{tikzpicture}
  \tikzstyle{var}=[circle,draw=black,fill=white,semithick,minimum size=12pt]
  \tikzstyle{fac}=[rectangle,draw=black,fill=white,semithick,minimum size=12pt]
  \small
  \node (a) at (-1,3.5) {(a)};
  \begin{scope}[scale=0.7]
  \node (psi) [fac] at (2,2) {$\Fac{}$};
  \node (D) [fac] at (2,-1.5) {?};
  \foreach \x in {1,2,3} {
    \node (A\x) [var] at (\x,4) {};
    \path (A\x) edge (psi);
    \node (B\x) [var] at (4,\x) {};
    \path (B\x) edge (psi);
    \node (C\x) [var] at (\x,0) {};
    \path (C\x) edge (psi);
    \path (C\x) edge (D);
  } 
  \draw[dotted] (2,4) ellipse (1.5cm and 0.75cm);
  \draw[xshift=-1.8cm] node () at (2,4) {$A$};
  \draw[dotted] (4,2) ellipse (0.75cm and 1.5cm);
  \draw[xshift=1cm] node () at (4,2) {$B$};
  \draw[dotted] (2,0) ellipse (1.5cm and 0.75cm);
  \draw[xshift=-1.8cm] node () at (2,0) {$C$};
  \end{scope}
  \node (b) at (6,3.5) {(b)};
  \begin{scope}[xshift=7cm,scale=0.7]
  \node (psi) [fac] at (2,2) {$\Fac{}$};
  \foreach \x in {1,2,3} {
    \node (A\x) [var] at (\x,4) {};
    \path (A\x) edge (psi);
    \node (B\x) [var] at (4,\x) {};
    \path (B\x) edge (psi);
    \node (C\x) [var] at (\x,0) {};
    \path (C\x) edge (psi);
    \node (D\x) [fac] at (\x,-1.5) {?};
    \path (C\x) edge (D\x);
  } 
  \draw[dotted] (2,4) ellipse (1.5cm and 0.75cm);
  \draw[xshift=-1.8cm] node () at (2,4) {$A$};
  \draw[dotted] (4,2) ellipse (0.75cm and 1.5cm);
  \draw[xshift=1cm] node () at (4,2) {$B$};
  \draw[dotted] (2,0) ellipse (1.5cm and 0.75cm);
  \draw[xshift=-1.8cm] node () at (2,0) {$C$};
  \end{scope}
\end{tikzpicture}
\caption{\label{fig:basic_lemma}The basic lemma: the smallest bounding
box enclosing the set of possible marginals of $\x{A}$ is identical in 
cases (a) and (b), if we are allowed to put arbitrary factors on the 
factor nodes marked with question marks.}
\end{figure}
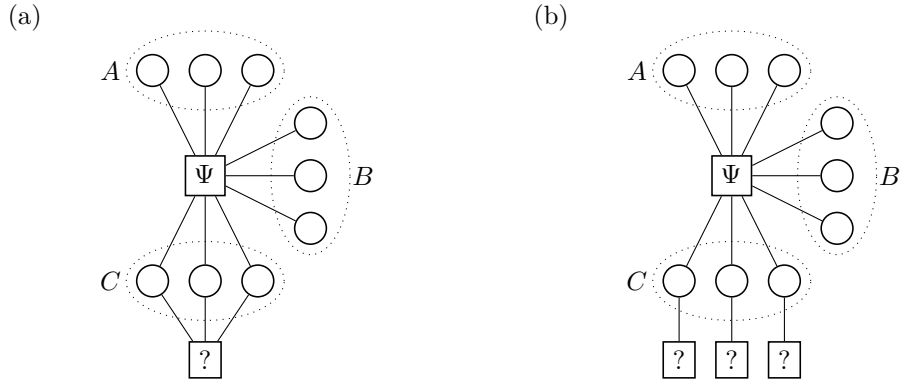

After defining the elementary concepts, we can proceed with the basic lemma.
Given the definitions introduced before, the basic lemma is easy to formulate.
It is illustrated in Figure \ref{fig:basic_lemma}.
\begin{lemm}\label{lemm:basic}
Let $A, B, C \subseteq \vars$ be mutually disjoint subsets of variables.
Let $\Fac{} \in \meas{A \cup B \cup C}$ such that for each $\x{C}\in\CX{C}$,
$$\sum_{\x{A\cup B}} \Fac{} > 0.$$ 
Then:
\begin{equation*}
\sbbox{\nor \left( \sum_{\x{B}, \x{C}} \Fac{} \meas{C}^* \right)}
= \sbbox{\nor \left( \sum_{\x{B}, \x{C}} \Fac{} \fmeas{C}^* \right)}.
\end{equation*}
\end{lemm}
\begin{proof}
Note that $\meas{C}^*$ is the convex hull of $\fmeas{C}^*$. Furthermore,
the multiplication with $\Fac{}$ and the summation over $\x{B}, \x{C}$
preserves convex combinations, as does the normalization operation (see Lemma
\ref{lemm:normalizing_preserves_convexity}). Therefore,
\begin{equation*}
\nor \left( \sum_{\x{B}, \x{C}} \Fac{} \meas{C}^* \right) \subseteq \hull \left( \nor \left( \sum_{\x{B}, \x{C}} \Fac{} \fmeas{C}^* \right) \right)
\end{equation*}
from which the lemma follows.
\end{proof}
The positivity condition is a technical condition, which in our experience is fulfilled for
most practically relevant factor graphs.

\subsection{Examples}

Before proceeding to the first main result, we first illustrate for a simple case
how the basic lemma
can be employed to obtain bounds on marginals. We show two bounds for the
marginal of the variable $\x{i}$ in the factor graph in Figure
\ref{fig:bound_example}(a).

\subsubsection{Example I}

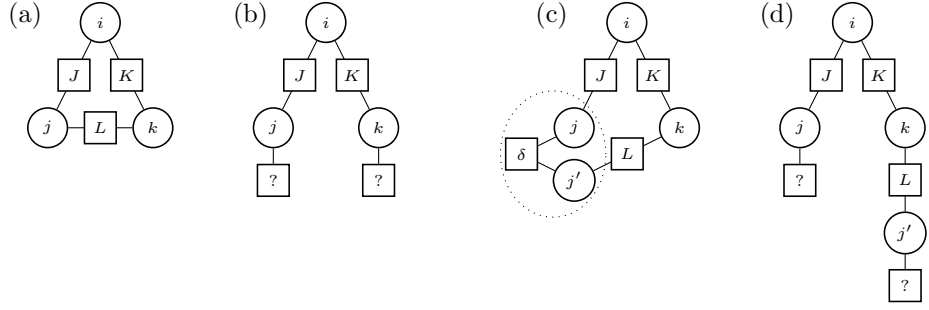
\begin{figure}
\centering
\begin{tikzpicture}
  \tikzstyle{var}=[circle,draw=black,fill=white,semithick,minimum size=15pt]
  \tikzstyle{fac}=[rectangle,draw=black,fill=white,semithick,minimum size=12pt]
  \small
  \node (a) at (-1,1.5) {(a)};
  \begin{scope}[xshift=0cm,scale=0.7]
    \tiny
    \node (i) at (0,2)  [var] {$i$};
    \node (j) at (-1,0) [var] {$j$};
    \node (k) at (1,0)  [var] {$k$};
    \path (i) edge node[pos=0.5] (J) [fac] {$J$} (j);
    \path (i) edge node[pos=0.5] (K) [fac] {$K$} (k);
    \path (j) edge node[pos=0.5] (L) [fac] {$L$} (k);
  \end{scope}
  \node (b) at (2,1.5) {(b)};
  \begin{scope}[xshift=3cm,scale=0.7]
    \tiny
    \node (i) at (0,2)  [var] {$i$};
    \node (j) at (-1,0) [var] {$j$};
    \node (k) at (1,0)  [var] {$k$};
    \path (i) edge node[pos=0.5] (J) [fac] {$J$} (j);
    \path (i) edge node[pos=0.5] (K) [fac] {$K$} (k);
    \node (jq) at (-1,-1) [fac] {?};
    \node (kq) at (1,-1) [fac] {?};
    \path (j) edge (jq);
    \path (k) edge (kq);
  \end{scope}
  \node (c) at (6,1.5) {(c)};
  \begin{scope}[xshift=7cm,scale=0.7]
    \tiny
    \node (i) at (0,2)  [var] {$i$};
    \node (j) at (-1,0) [var] {$j$};
    \node (k) at (1,0)  [var] {$k$};
    \path (i) edge node[pos=0.5] (J) [fac] {$J$} (j);
    \path (i) edge node[pos=0.5] (K) [fac] {$K$} (k);
    \node (delta) at (-2,-0.5) [fac] {$\delta$};
    \node (jq) at (-1,-1) [var] {$j'$};
    \path (j)  edge (delta);
    \path (jq) edge (delta);
    \path (k) edge node[pos=0.5] (L) [fac] {$L$} (jq);
    \draw [dotted] (-1.4,-0.5) ellipse (1.0cm and 1.2cm);
  \end{scope}
  \node (d) at (9,1.5) {(d)};
  \begin{scope}[xshift=10cm,scale=0.7]
    \tiny
    \node (i) at (0,2)  [var] {$i$};
    \node (j) at (-1,0) [var] {$j$};
    \node (k) at (1,0)  [var] {$k$};
    \path (i) edge node[pos=0.5] (J) [fac] {$J$} (j);
    \path (i) edge node[pos=0.5] (K) [fac] {$K$} (k);
    \node (jq) at (-1,-1) [fac] {?};
    \node (kk) at (1,-2) [var] {$j'$};
    \path (k) edge node[pos=0.5] (L) [fac] {$L$} (kk);
    \node (kq) at (1,-3) [fac] {?};
    \path (kk) edge (kq);
    \path (j) edge (jq);
  \end{scope}
\end{tikzpicture}
\caption{\label{fig:bound_example}(a) Factor graph; 
(b) Illustration of the bound on $\Prob(\x{i})$ corresponding to Example I; 
(c) Cloning node $j$ by adding a new variable $j'$ and a factor $\fac{\delta}(\x{j},\x{j'}) = \delta_{\x{j}}(\x{j'})$; 
(d) Illustration of the improved bound on $\Prob(\x{i})$, corresponding to Example (II), based on (c).}
\end{figure}

First, note that the marginal of $\x{i}$ satisfies
\begin{equation*}
\Prob(\x{i}) 
= \nor \left( \sum_{\x{j}} \sum_{\x{k}} \fac{J} \fac{K} \fac{L} \right) 
\in \nor \left( \sum_{\x{j}} \sum_{\x{k}} \fac{J} \fac{K} \meas{\{j,k\}}^* \right).
\end{equation*}
because, obviously, $\fac{L} \in \meas{\{j,k\}}^*$. Now, applying the basic lemma
with $A = \{i\}$, $B = \emptyset$, $C = \{j,k\}$ and $\Fac{} = \fac{J}\fac{K}$, we obtain
\begin{equation*}
\Prob(\x{i}) 
\in \sbbox{\nor \left( \sum_{\x{j}} \sum_{\x{k}} \fac{J} \fac{K} \fmeas{\{j,k\}}^* \right)}.
\end{equation*}
Applying the distributive law, we conclude
\begin{equation*}
\Prob(\x{i}) 
\in \snbbox{\Big( \sum_{\x{j}} \fac{J} \meas{j}^* \Big) \Big( \sum_{\x{k}} \fac{K} \meas{k}^* \Big)},
\end{equation*}
which certainly implies
\begin{equation*}
\Prob(\x{i}) 
\in \snbbox{\snbbox{\sum_{\x{j}} \fac{J} \meas{j}^*} \cdot \snbbox{\sum_{\x{k}} \fac{K} \meas{k}^*}}.
\end{equation*}
This is illustrated in Figure \ref{fig:bound_example}(b), which should be read as ``What can we
say about the range of $\Prob(\x{i})$ when the factors corresponding to the nodes marked with 
question marks are arbitrary?''
Because of the various occurrences of the normalization operator, we can restrict ourselves
to normalized measures on the question-marked factor nodes:
\begin{equation*}
\Prob(\x{i}) 
\in \snbbox{\snbbox{\sum_{\x{j}} \fac{J} \nmeas{j}} \cdot \snbbox{\sum_{\x{k}} \fac{K} \nmeas{k}}}.
\end{equation*}
Now it may seem that this smallest bounding box would be difficult to compute, because in
principle one would have to compute all the measures in the sets $\nor\sum_{\x{j}} \fac{J} \nmeas{j}$
and $\nor\sum_{\x{k}} \fac{K} \nmeas{k}$. Fortunately, we only need to compute the extreme points of 
these sets, because the mapping 
$$\meas{\{j\}}^* \to \meas{\{i\}}^* : \fac{} \mapsto \nor \sum_{\x{j}} \fac{J} \fac{}$$
maps convex combinations into convex combinations (and similarly for the other mapping, involving
$\fac{K}$). Since smallest bounding boxes only depend on extreme points, we conclude that
\begin{equation*}
\Prob(\x{i}) 
\in \snbbox{\snbbox{\sum_{\x{j}} \fac{J} \ep \nmeas{j}} \cdot \snbbox{\sum_{\x{k}} \fac{K} \ep \nmeas{k}}}
\end{equation*}
which can be calculated efficiently if the number of possible values of each variable is small.

\subsubsection{Example II}

We can improve this bound by using another trick: cloning variables.
The idea is to first clone the variable $\x{j}$ by adding a new variable $\x{j'}$ that is
constrained to take the same value as $\x{j}$. In terms of the factor graph, we add a
variable node $j'$ and a factor node $\delta$, connected to variable nodes $j$ and $j'$, 
with corresponding factor $\fac{\delta}(\x{j},\x{j'}) := \delta_{\x{j}}(\x{j'})$;
see also Figure \ref{fig:bound_example}(c). Clearly, the marginal of $\x{i}$ satisfies:
\begin{equation*}\begin{split}
\Prob(\x{i}) 
& = \nor \left( \sum_{\x{j}} \sum_{\x{k}} \fac{J} \fac{K} \fac{L} \right) \\
& = \nor \left( \sum_{\x{j}} \sum_{\x{j'}} \sum_{\x{k}} \fac{J} \fac{K} \fac{L} \delta_{\x{j}}(\x{j'}) \right) \\
\end{split}\end{equation*}
where it should be noted that in the first line, $\fac{L}$ is shorthand for
$\fac{L}(\x{j},\x{k})$ but in the second line it is meant as shorthand for
$\fac{L}(\x{j'},\x{k})$. Noting that $\fac{\delta} \in \meas{\{j,j'\}}^*$ and applying the
basic lemma with $C = \{j,j'\}$ yields:
\begin{equation*}
\Prob(\x{i}) 
 \in \nor \left( \sum_{\x{j}} \sum_{\x{j'}} \sum_{\x{k}} \fac{J} \fac{K} \fac{L} \meas{\{j,j'\}}^* \right) 
 \in \snbbox{\sum_{\x{j}} \sum_{\x{j'}} \sum_{\x{k}} \fac{J} \fac{K} \fac{L} \fmeas{\{j,j'\}}^*}.
\end{equation*}
Applying the distributive law, we obtain (see also Figure \ref{fig:bound_example}(d)):
\begin{equation*}
\Prob(\x{i}) 
\in \snbbox{\left(\sum_{\x{j}} \fac{J} \meas{\{j\}}^* \right) \left( \sum_{\x{k}} \fac{K} \sum_{\x{j'}} \fac{L} \meas{\{j'\}}^* \right)},
\end{equation*}
from which we conclude
\begin{equation*}
\Prob(\x{i}) 
\in \snbbox{\snbbox{\sum_{\x{j}} \fac{J} \nmeas{\{j\}}} \snbbox{\sum_{\x{k}} \fac{K} \snbbox{\sum_{\x{j'}} \fac{L} \nmeas{\{j'\}}}}}.
\end{equation*}
This can again be calculated efficiently by considering only extreme points.

As a more concrete example, take all variables as binary and take for each
factor $\fac{} = \big( \begin{smallmatrix} 1 & 2 \\ 2 & 1 \end{smallmatrix} \big)$.
Then the first bound (Example I) yields:
$$\Prob(\x{i}) \in \bbox{i}{\Big(\begin{smallmatrix}1/5 \\ 1/5\end{smallmatrix}\Big)}{\Big(\begin{smallmatrix}4/5 \\ 4/5\end{smallmatrix}\Big)},$$
whereas the second, tighter, bound (Example II) gives:
$$\Prob(\x{i}) \in \bbox{i}{\Big(\begin{smallmatrix}2/7 \\ 2/7\end{smallmatrix}\Big)}{\Big(\begin{smallmatrix}5/7 \\ 5/7\end{smallmatrix}\Big)}.$$
Obviously, the exact marginal is 
$$\Prob(\x{i}) = \Big(\begin{smallmatrix}1/2 \\ 1/2\end{smallmatrix}\Big).$$


\subsection{Propagation of boxes over a subtree}

\begin{figure}
\centering
\begin{tikzpicture}
  \tikzstyle{var}=[circle,draw=black,fill=white,semithick,minimum size=0.5cm]
  \tikzstyle{fac}=[rectangle,draw=black,fill=white,semithick,minimum size=0.4cm]
  \tikzstyle{arr}=[->,shorten <=1pt,>=stealth']
  \small
  \node (a) at (-1,2) {(a)};
  \begin{scope}[xshift=0.5cm,scale=0.7]
    \tiny
    \node (i) at (0,2)  [var] {$i$};
    \node (j) at (-1,0) [var] {$j$};
    \node (k) at (1,0)  [var] {$k$};
    \path (i) edge node[pos=0.5] (J) [fac] {$J$} (j);
    \path (i) edge node[pos=0.5] (K) [fac] {$K$} (k);
    \path (j) edge node[pos=0.5] (L) [fac] {$L$} (k);
  \end{scope}
  \node (b) at (3,2) {(b)};
  \begin{scope}[xshift=4.5cm,yshift=0cm,scale=0.7]
    \tiny
    \node (i) at (0,2)  [var] {$i$};
    \node (j) at (-1,0) [var] {$j$};
    \node (k) at (1,0)  [var] {$k$};
    \path (i) edge node[pos=0.5] (J) [fac] {$J$} (j);
    \path (i) edge node[pos=0.5] (K) [fac] {$K$} (k);
    \node (L) at (1,-1) [fac] {L};
    \path (k) edge (L);
  \end{scope}
  \node (c) at (7,2) {(c)};
  \begin{scope}[xshift=9cm,yshift=0cm,scale=0.7]
    \tiny
    \node (i) at (0,2)  [var] {$i$};
    \node (j) at (-1.6,-1.2) [var] {$j$};
    \node (k) at (1.6,-1.2)  [var] {$k$};
    \path (i) edge node[pos=0.5] (J) [fac] {$J$} (j);
    \path (i) edge node[pos=0.5] (K) [fac] {$K$} (k);
    \node (L) at (1.6,-2.8) [fac] {L};
    \path (k) edge (L);
    \node (LL) at (-1.6,-2.8) {$\mathcal{P}_j$};
    \node (jj) at (1.6,-4.4) {$\mathcal{P}_j$};
    \draw [arr] (J) -- (i);
    \draw [arr] (K) -- (i);
    \draw [arr] (j) -- (J);
    \draw [arr] (k) -- (K);
    \draw [arr] (L) -- (k);
    \node (Bi)  at (0,2.6) {$\mathcal{B}_i$};
    \node (BJi) at (-1.0,1.2) {$\mathcal{B}_{J \to i}$};
    \node (BKi) at (1.0,1.2) {$\mathcal{B}_{K \to i}$};
    \node (BjJ) at (-1.8,-0.4) {$\mathcal{B}_{j \to J}$};
    \node (BkK) at (1.85,-0.4) {$\mathcal{B}_{k \to K}$};
    \node (BLk) at (2.2,-2.1) {$\mathcal{B}_{L \to k}$};
    \node (BLj) at (-2.2,-2.1) {$\mathcal{B}_{L \to j}$};
    \node (BjL) at (2.2,-3.7) {$\mathcal{B}_{j \to L}$};
    \draw [arr] (-1.6,-2.55) -- (j);
    \draw [arr] (1.6,-4.15) -- (L);
  \end{scope}
\end{tikzpicture}
\caption{\label{fig:boxprop_subtree}Box propagation algorithm corresponding to Example II: 
(a) Factor graph $\graph = (\vars,\factors,\uedges)$; 
(b) a possible subtree $(V,F,E)$ of $G$; 
(c) propagating sets of measures (boxes or simplices) on the subtree leading to a bound $\bebox{i}$ 
on the marginal probability of $\x{i}$ in $\graph$.}
\end{figure}
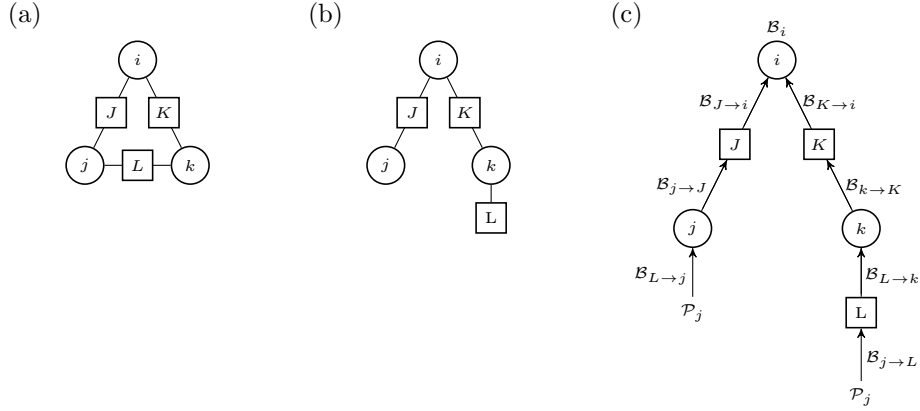

We now formulate a message passing algorithm that
resembles Belief Propagation. However, instead of propagating measures, it
propagates boxes (or simplices) of measures; furthermore, it is only applied to a
subtree of the factor graph, propagating boxes from the leaves towards a root
node, instead of propagating iteratively over the whole factor graph several
times. The resulting ``belief'' at the root node is a box that bounds the
exact marginal of the root node. The choice of the subtree is arbitrary;
different choices lead to different bounds in general. We illustrate the 
algorithm using the example that we have studied before (see Figure \ref{fig:boxprop_subtree}).

\begin{defi}
Let $(\vars, \factors, \uedges)$ be a factor graph. We call the bipartite
graph $(V, F, E)$ a \emph{subtree of $(\vars,\factors,\uedges)$ with root $i$}
if $i \in V \subseteq \vars$, $F \subseteq \factors$, $E \subseteq \uedges$
such that $(V,F,E)$ is a tree with root $i$ and for all 
$\ue{j}{J} \in E$, $j \in V$ and $J \in F$ (\ie there are no ``loose edges'').\footnote{Note 
that this corresponds to the notion of subtree of a bipartite graph;
for a subtree of a factor graph, one sometimes imposes the additional 
constraint that for all factors $J \in F$, all its connecting edges 
$\ue{J}{j}$ with $j \in \nbf{J}$
have to be in $E$; here we do not impose this additional constraint.}
\end{defi}

An illustration of a factor graph and a possible subtree is given in Figure \ref{fig:boxprop_subtree}(a)-(b).
We denote the parent of $j \in V$ according to $(V,F,E)$ by $\parents{j}$ and
similarly, we denote the parent of $J \in F$ by $\parents{J}$. In the following, we will use
the topology of the \emph{original} factor graph $(\vars,\factors,\uedges)$ whenever we refer to
neighbors of variables or factors.

Each edge of the subtree will carry one message, oriented such that it
``flows'' towards the root node. In addition, we define messages entering the
subtree for all ``missing'' edges in the subtree.  Because of the
bipartite character of the factor graph, we can distinguish between two types
of messages: messages $\mebox{J}{j} \subseteq \meas{j}$ 
sent to a variable $j \in V$ from a neighboring factor $J \in \nbv{j}$, and messages
$\mebox{j}{J}\subseteq \meas{j}$ sent to a factor $J \in F$ from a neighboring
variable $j\in\nbf{J}$. 

The messages entering the subtree are all defined to be simplices; more precisely,
we define the incoming messages
\begin{align*}
\mebox{j}{J} & = \nmeas{j} \qquad \text{$J \in F$, $\ue{j}{J} \in \uedges\setm E$}\\
\mebox{J}{j} & = \nmeas{j} \qquad \text{$j \in V$, $\ue{j}{J} \in \uedges\setm E$}.
\end{align*}
We propagate messages towards the root $i$ of the tree using the following 
update rules (note the similarity with the BP update rules). The message sent
from a variable $j \in V$ to its parent $J = \parents{j} \in F$ is defined as
\begin{equation*}
\mebox{j}{J} = \begin{cases}
\displaystyle\prod_{K \in \nbve{j}{J}} \mebox{K}{j} & \text {if all incoming $\mebox{K}{j}$ are boxes} \\[5pt]
\nmeas{j} & \text{if at least one of the $\mebox{K}{j}$ is the simplex $\nmeas{j}$},
\end{cases}
\end{equation*}
where the product of the boxes can be calculated using Lemma \ref{lemm:product_boxes_same_varset}.
The message sent from a factor $J \in F$ to its parent $k = \parents{J} \in V$ is defined as
\begin{equation}\label{eq:boxprop_subtree_factor_to_variable}
\mebox{J}{k} = \snbbox{\sum_{\x{\nbfe{J}{k}}} \fac{J} \prod_{l \in \nbfe{J}{k}} \mebox{l}{J}}.
\end{equation}
This smallest bounding box can be calculated using
the following Corollary of Lemma \ref{lemm:product_convex_sets}:
\begin{coro}
\begin{equation*}
\snbbox{\sum_{\x{\nbfe{J}{k}}} \fac{J} \prod_{l \in \nbfe{J}{k}} \mebox{l}{J}} =
\snbbox{\sum_{\x{\nbfe{J}{k}}} \fac{J} \prod_{l \in \nbfe{J}{k}} \ep \mebox{l}{J}}
\end{equation*}
\end{coro}
\begin{proof}
By Lemma \ref{lemm:product_convex_sets}, 
\begin{equation*}
\prod_{l \in \nbfe{J}{k}} \mebox{l}{J} \subseteq \hull \left( \prod_{l \in \nbfe{J}{k}} \ep \mebox{l}{J} \right).
\end{equation*}
Because the multiplication with $\fac{J}$ and the summation over
$\x{\nbfe{J}{k}}$ preserves convex combinations, as does the normalization
(see Lemma \ref{lemm:normalizing_preserves_convexity}),
the statement follows.
\end{proof}
The final ``belief'' $\bebox{i}$
at the root node $i$ is calculated by
\begin{equation*}
\bebox{i} = \begin{cases}
\displaystyle\snbbox{\prod_{K \in \nbv{j}} \mebox{K}{j}} & \text {if all incoming $\mebox{K}{j}$ are boxes} \\[5pt]
\nmeas{j} & \text{if at least one of the $\mebox{K}{j}$ is the simplex $\nmeas{j}$}.
\end{cases}
\end{equation*}
Note that when applying this to the case illustrated in Figure \ref{fig:boxprop_subtree}, we
obtain the bound that we derived earlier on (``Example II'').

We can now formulate our first main result, which gives a rigorous bound on the
exact single-variable marginal of the root node:
\begin{theo}\label{theo:main}
Let $(\vars,\factors,\uedges)$ be a factor graph with corresponding 
probability distribution \eref{eq:prob_dist}.
Let $i \in \vars$ and $(V,F,E)$ be a subtree of $(\vars,\factors,\uedges)$ 
with root $i \in V$. Apply the ``box propagation'' algorithm described
above to calculate the final ``belief'' $\bebox{i}$ on the root node
$i$. Then $\Prob(\x{i}) \in \bebox{i}$.
\end{theo}
\begin{proofsketch}
The first step consists in 
extending the subtree such that each factor node has the right number of
neighboring variables by cloning the missing variables. The second step consists 
of applying the basic lemma where the set $C$ consists of all the variable 
nodes of the subtree which have connecting edges in $\uedges\setm E$, together
with all the cloned variable nodes. Then we apply the distributive law, 
which can be done because the extended subtree has no cycles. Finally, we relax
the bound by adding additional normalizations and smallest bounding boxes at each
factor node in the subtree. It should now be clear that the recursive algorithm 
``box propagation'' described above precisely calculates the smallest bounding 
box at the root node $i$ that corresponds to this procedure.
\end{proofsketch}

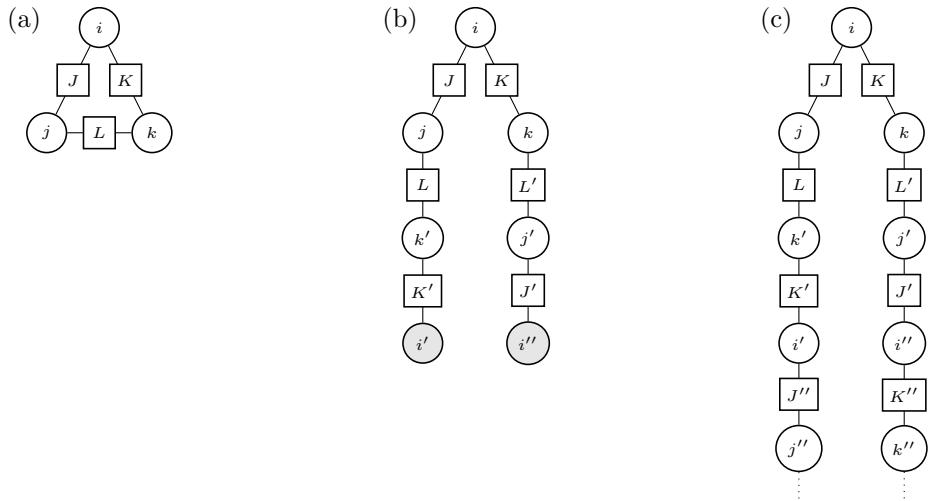
\begin{figure}
\centering
\begin{tikzpicture}
  \tikzstyle{var}=[circle,draw=black,fill=white,semithick,minimum size=15pt]
  \tikzstyle{fac}=[rectangle,draw=black,fill=white,semithick,minimum size=12pt]
  \small
  \node (a) at (-1,1.5) {(a)};
  \begin{scope}[xshift=0cm,scale=0.7]
    \tiny
    \node (i) at (0,2)  [var] {$i$};
    \node (j) at (-1,0) [var] {$j$};
    \node (k) at (1,0)  [var] {$k$};
    \path (i) edge node[pos=0.5] (J) [fac] {$J$} (j);
    \path (i) edge node[pos=0.5] (K) [fac] {$K$} (k);
    \path (j) edge node[pos=0.5] (L) [fac] {$L$} (k);
  \end{scope}
  \node (b) at (4,1.5) {(b)};
  \begin{scope}[xshift=5cm,scale=0.7]
    \tiny
    \node (i)  at (0,2)   [var] {$i$};
    \node (j)  at (-1,0)  [var] {$j$};
    \node (k)  at (1,0)   [var] {$k$};
    \path (i)  edge node[pos=0.5] (J) [fac] {$J$} (j);
    \path (i)  edge node[pos=0.5] (K) [fac] {$K$} (k);
    \node (k2) at (-1,-2) [var] {$k'$};
    \node (j2) at (1,-2)  [var] {$j'$};
    \path (j)  edge node[pos=0.5] (L)  [fac] {$L$}  (k2);
    \path (k)  edge node[pos=0.5] (L2) [fac] {$L'$} (j2);
    \node (i2) at (-1,-4) [var,fill=gray!20] {$i'$};
    \node (i3) at (1,-4)  [var,fill=gray!20] {$i''$};
    \path (k2) edge node[pos=0.5] (K2) [fac] {$K'$} (i2);
    \path (j2) edge node[pos=0.5] (J2) [fac] {$J'$} (i3);
  \end{scope}
  \node (c) at (9,1.5) {(c)};
  \begin{scope}[xshift=10cm,scale=0.7]
    \tiny
    \node (i)  at (0,2)   [var] {$i$};
    \node (j)  at (-1,0)  [var] {$j$};
    \node (k)  at (1,0)   [var] {$k$};
    \path (i)  edge node[pos=0.5] (J) [fac] {$J$} (j);
    \path (i)  edge node[pos=0.5] (K) [fac] {$K$} (k);
    \node (k2) at (-1,-2) [var] {$k'$};
    \node (j2) at (1,-2)  [var] {$j'$};
    \path (j)  edge node[pos=0.5] (L)  [fac] {$L$}  (k2);
    \path (k)  edge node[pos=0.5] (L2) [fac] {$L'$} (j2);
    \node (i2) at (-1,-4) [var] {$i'$};
    \node (i3) at (1,-4)  [var] {$i''$};
    \path (k2) edge node[pos=0.5] (K2) [fac] {$K'$} (i2);
    \path (j2) edge node[pos=0.5] (J2) [fac] {$J'$} (i3);
    \node (j3) at (-1,-6) [var] {$j''$};
    \node (k3) at (1,-6)  [var] {$k''$};
    \path (i2) edge node[pos=0.5] (K3) [fac] {$J''$} (j3);
    \path (i3) edge node[pos=0.5] (J3) [fac] {$K''$} (k3);
    \path (j3) edge [dotted] (-1,-7);
    \path (k3) edge [dotted] (1,-7);
  \end{scope}
\end{tikzpicture}
\caption{\label{fig:trees_example}(a) Factor graph; 
(b) Self-avoiding walk tree with root $i$, with cycle-induced leaf nodes shown in gray; 
(c) Computation tree for $i$.}
\end{figure}

Note that a subtree of the orginal factor graph is also a subtree of the
\emph{computation tree for $i$} \citep{TatikondaJordan02}. A computation tree
is an ``unwrapping'' of the factor graph that has been used in analyses of the
Belief Propagation algorithm. The computation tree starting at variable $i \in
\vars$ consists of all paths on the factor graph, starting at $i$, that never
backtrack (see also Figure \ref{fig:trees_example}(c)). This means that the
bounds on the (exact) marginals that we just derived are at the same time
bounds on the approximate Belief Propagation marginals (beliefs).

\begin{coro}
In the situation described in Theorem \ref{theo:main}, the final bounding 
box $\bebox{i}$ also bounds the (approximate) Belief Propagation 
marginal of the root node $i$, \ie $\Prob_{BP}(\x{i}) \in \bebox{i}$.\qed
\end{coro}

\subsection{Bounds using Self-Avoiding Walk Trees}

While writing this article, we became aware that a related method to
obtain bounds on single-variable marginals has been proposed recently by
\citet{Ihler07}.\footnote{Note that \citep[Lemma 5]{Ihler07} contains an error:
to obtain the correct expression, one has to replace $\delta$ with $\delta^2$,
\ie the correct statement would be that 
$$\frac{m(j)}{\delta^2 + (1-\delta^2)m(j)} \le p(j) \le \frac{\delta^2 m(j)}{1 - (1-\delta^2) m(j)}$$ 
if $d\big(p(x)/m(x)\big) \le \delta$ (where $p$ and $m$ should both be
normalized).} The method presented there uses a different local bound, which
empirically seems to be less tight than ours, but has the advantage of 
being computationally less demanding if the domains of the random variables are large.
On the other hand, the
bound presented there does not use subtrees of the factor graph, but 
uses self-avoiding walk (SAW) trees instead. Since each subtree of the factor graph
is a subtree of an SAW tree, this may lead to tighter bounds.

The idea of using a self-avoiding walk tree for calculating marginal
probabilities seems to be generally attributed to \citet{Weitz06}, but
can already be found in \citep{ScottSokal05}. In this
subsection, we show how this idea can be combined with the propagation of
bounding boxes. The result Theorem \ref{theo:saw} will turn out to be an improvement over Theorem
\ref{theo:main} in case there are only pairwise interactions, whereas in the
general case, Theorem \ref{theo:main} often yields tighter bounds
empirically.

\begin{defi}
Let $\graph := (\vars, \factors, \uedges)$ be a factor graph and let $i \in \vars$. 
A \emph{self-avoiding walk (SAW)} starting at $i \in \vars$ of length $n \in \NN^*$
is a sequence $\alpha = (\alpha_1, \alpha_2, \alpha_3, \dots, \alpha_n) \in (\vars \cup \factors)^n$ that
{\renewcommand{\labelenumi}{(\roman{enumi})}
\begin{enumerate}
\item starts at $i \in \vars$, \ie $\alpha_1 = i$;
\item subsequently visits neighboring nodes in the factor graph, \ie $\alpha_{j+1} \in \nb{\alpha_j}$ for all $j = 1, 2, \dots, n-1$;
\item does not backtrack, \ie $\alpha_j \ne \alpha_{j+2}$ for all $j = 1, 2, \dots, n-2$;
\item the first $n-1$ nodes are all different, \ie
$\alpha_j \ne \alpha_k$ if $j\ne k$ for $j,k \in \{1,2,\dots,n-1\}$.\footnote{Note that (iii) almost follows from (iv), except for
the condition that $\alpha_{n-2} \ne \alpha_n$.}
\end{enumerate}}
The set of all self-avoiding walks starting at $i \in \vars$ has a natural tree structure,
defined by declaring each SAW $(\alpha_1, \alpha_2, \dots, \alpha_n, \alpha_{n+1})$ to be 
a child of the SAW $(\alpha_1, \alpha_2, \dots, \alpha_n)$, for all $n \in \NN^*$; the 
resulting tree is called the \emph{self-avoiding walk (SAW) tree with root $i \in \vars$},
denoted $T^{SAW}_{\graph}(i)$.
\end{defi}

Note that the name ``self-avoiding walk tree'' is slightly inaccurate, because
the last node of a SAW may have been visited already. In general, the SAW tree
can be much larger than the original factor graph.  Following \citet{Ihler07},
we call a leaf node in the SAW tree a \emph{cycle-induced leaf node} if it
contains a cycle (\ie if its final node has been visited before in the same walk), and call it a 
\emph{dead-end leaf node} otherwise.  We denote the parent of node $\alpha$ in the SAW tree
by $\parents{\alpha}$ and we denote its children by $\children{\alpha}$.
The final node of a SAW $\alpha = (\alpha_1,\dots,\alpha_n)$ is denoted by 
$\graph(\alpha) = \alpha_n$. An example of a SAW tree for our running example
factor graph is shown in Figure \ref{fig:trees_example}(b).

Let $\graph = (\vars,\factors,\uedges)$ be a factor graph and let $i \in \vars$.
We now define a 
propagation algorithm on the SAW tree $T^{SAW}_{\graph}(i)$, where each node 
$\alpha \in T^{SAW}_{\graph}(i)$ (except for the root $i$) sends a message 
$\mebox{\alpha}{\parents{\alpha}}$ to 
its parent node $\parents{\alpha} \in T^{SAW}_{\graph}(i)$. Each cycle-induced 
leaf node of $T^{SAW}_{\graph}(i)$ sends a simplex
to its parent node: if $\alpha$ is a cycle-induced leaf node, then
\begin{equation}\label{eq:boxprop_saw_cycle_induced_leaf_node}
\mebox{\alpha}{\parents{\alpha}} = \begin{cases}
\nmeas{\graph(\alpha)} & \text{if $\graph(\alpha) \in \vars$} \\[5pt]
\nmeas{\graph(\parents{\alpha})} & \text{if $\graph(\alpha) \in \factors$.}
\end{cases}
\end{equation}
All other nodes $\alpha$ in the SAW tree (\ie the dead-end leaf nodes and
the nodes with children, except for the root $i$) send a message according 
to the following rules. If $\graph(\alpha) \in \vars$, 
\begin{equation}\label{eq:boxprop_saw_variable_to_factor}
\mebox{\alpha}{\parents{\alpha}} = \begin{cases}
\displaystyle\prod_{\beta \in \children{\alpha}} \mebox{\beta}{\alpha} & \text {if all $\mebox{\beta}{\alpha}$ are boxes} \\[5pt]
\nmeas{\graph(\alpha)} & \text{if at least one of the $\mebox{\beta}{\alpha}$ is a simplex.}
\end{cases}
\end{equation}
On the other hand, if $\graph(\alpha) \in \factors$,
\begin{equation}\label{eq:boxprop_saw_factor_to_variable}
\mebox{\alpha}{\parents{\alpha}} = \snbbox{\sum_{\x{\setm \graph(\parents{\alpha})}} \fac{\graph(\alpha)} \sbbox{\prod_{\beta \in \children{\alpha}} \mebox{\beta}{\alpha}}}.
\end{equation}
The final ``belief'' at the root node $i \in \vars$ is defined as:
\begin{equation}\label{eq:boxprop_saw_final_result}
\bebox{i} = \begin{cases}
\displaystyle\snbbox{\prod_{\beta \in \children{i}} \mebox{\beta}{i}} & \text {if all $\mebox{\beta}{i}$ are boxes} \\[5pt]
\nmeas{\graph(i)} & \text{if at least one of the $\mebox{\beta}{i}$ is a simplex.}
\end{cases}
\end{equation}
We will refer to this algorithm as ``box propagation on the SAW tree''; 
it is similar to the propagation algorithm for boxes on subtrees of the factor
graph that we defined earlier. However, note that 
whereas \eref{eq:boxprop_subtree_factor_to_variable} bounds a sum-product assuming
that incoming measures factorize, \eref{eq:boxprop_saw_factor_to_variable} is a looser bound that also
holds if the incoming measures do not necessarily factorize. In the special case where the factor
depends only on two variables, the updates \eref{eq:boxprop_subtree_factor_to_variable} and \eref{eq:boxprop_saw_factor_to_variable} are identical.

\begin{theo}\label{theo:saw}
Let $\graph := (\vars, \factors, \uedges)$ be a factor graph.
Let $i \in \vars$ and let $T^{SAW}_{\graph}(i)$ be the
SAW tree with root $i$. Then $\Prob(\x{i}) \in \bebox{i}$,
where $\bebox{i}$ is the bounding box that results from propagating bounds 
on the SAW tree $T^{SAW}_{\graph}(i)$ according to equations \eref{eq:boxprop_saw_cycle_induced_leaf_node}--\eref{eq:boxprop_saw_final_result}.
\end{theo}

\begin{figure}
\centering
\begin{tikzpicture}
  \tikzstyle{var}=[circle,draw=black,fill=white,semithick,minimum size=12pt]
  \tikzstyle{fac}=[rectangle,draw=black,fill=white,semithick,minimum size=12pt]
  \small
  \node (a) at (-1,4.5) {(a)};
  \begin{scope}[scale=0.7]
  \node (psi) [fac] at (2.4,2) {$\Psi$};
  \tiny
  \node (i1) [var] at (0,4) {$i_1$};
  \node (i2) [var] at (1.2,4) {$i_2$};
  \node (i3) [var] at (2.4,4) {$i_3$};
  \node (ii) at (3.6,4) {$\dots$};
  \node (in) [var] at (4.8,4) {$i_n$};
  \path (i1) edge (psi);
  \path (i2) edge (psi);
  \path (i3) edge (psi);
  \path (in) edge (psi);
  \node (j1) [var] at (0,0) {};
  \node (j2) [var] at (1.6,0) {};
  \node (j3) at (3.2,0) {$\dots$};
  \node (jm) [var] at (4.8,0) {};
  \path (j1) edge (psi);
  \path (j2) edge (psi);
  \path (jm) edge (psi);
  \draw[dotted] (2.4,4) ellipse (3.4cm and 0.8cm);
  \draw[xshift=-4.3cm] node () at (3,4) {$A$};
  \draw[dotted] (2.4,0) ellipse (3.4cm and 0.8cm);
  \draw[xshift=-4.3cm] node () at (3,0) {$C$};
  \end{scope}
  \node (b) at (5.5,4.5) {(b)};
  \begin{scope}[xshift=6.5cm,scale=0.7]
  \node (psi) [fac] at (2.4,2) {$\Psi$};
  \tiny
  \node (i1) [var] at (0,4) {$i_1$};
  \node (i2) [var] at (1.2,4) {$i_2$};
  \node (i3) [var] at (2.4,4) {$i_3$};
  \node (ii) at (3.6,4) {$\dots$};
  \node (in) [var] at (4.8,4) {$i_n$};
  \path (i1) edge (psi);
  \path (i2) edge (psi);
  \path (i3) edge (psi);
  \path (in) edge (psi);
  \node (j1) [var] at (0,0) {};
  \node (j2) [var] at (1.6,0) {};
  \node (j3) at (3.2,0) {$\dots$};
  \node (jm) [var] at (4.8,0) {};
  \path (j1) edge (psi);
  \path (j2) edge (psi);
  \path (jm) edge (psi);
  \draw[dotted] (2.4,4) ellipse (3.4cm and 0.8cm);
  \draw[xshift=-4.3cm] node () at (3,4) {$A$};
  \draw[dotted] (2.4,0) ellipse (3.4cm and 0.8cm);
  \draw[xshift=-4.3cm] node () at (3,0) {$C$};
  \node (I2) [fac] at (1.2,5.5) {?};
  \node (I3) [fac] at (2.4,5.5) {?};
  \node (In) [fac] at (4.8,5.5) {?};
  \path (I2) edge (i2);
  \path (I3) edge (i3);
  \path (In) edge (in);
  \end{scope}
\end{tikzpicture}
\caption{\label{fig:second_basic_lemma}The second basic lemma: the 
marginal on $\x{A}$ in (a) is contained in the bounding box of the product of smallest
bounding boxes $B_{i}$ for $i \in A$, where (b) the smallest bounding box 
$B_i$ is obtained by putting arbitrary factors on the other variables 
in $A\setm \{i\}$ and calculating the smallest bounding box on $i$, illustrated here for the case $i=i_1$.}
\end{figure}
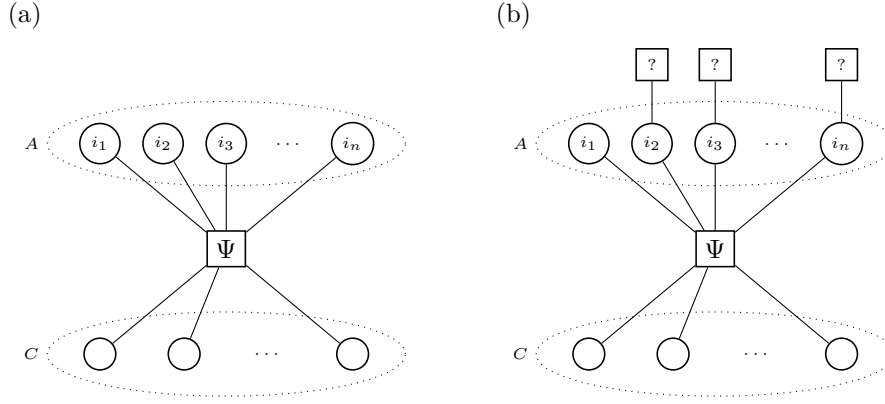

The following lemma, illustrated in Figure \ref{fig:second_basic_lemma},
plays a crucial role in the proof of the theorem. It seems to be related
to the so-called ``telegraph expansion'' used in \cite{Weitz06}.

\begin{lemm}
Let $A, C \subseteq \vars$ be two disjoint sets of variable indices and let 
$\Psi \in \meas{A \cup C}$ be a factor 
depending on (some of) the variables in $A \cup C$.
Then:
\begin{equation*}
\nor \left(\sum_{\x{C}} \Psi\right) \in \sbbox{\prod_{i\in A} B_i}
\end{equation*}
where
$$B_i := \snbbox{\sum_{\x{A\setm i}} \sum_{\x{C}} \Psi \fmeas{A\setm i}^*}.$$
\end{lemm}
\begin{proof}
We assume that $C = \emptyset$;
the more general case then follows from this special case by replacing $\Psi$ by
$\sum_{\x{C}} \Psi$. 

Let $A = \{i_1,i_2,\dots,i_n\}$ and let $\lFac{i}$, $\uFac{i}$ be the lower and
upper bounds corresponding to $B_{i}$,
for all $i\in A$. For each $k=1,2,\dots,n$, note that
$$\left(\prod_{l=1}^{k-1} \mathbf{1}_{i_l}\right) \left(\prod_{l=k+1}^n \delta_{\x{i_l}}\right) \in \fmeas{A\setm i_k}^*,$$
for all $\x{\{i_{k+1},\dots,i_n\}} \in \CX{\{i_{k+1},\dots,i_n\}}$. Therefore, we obtain from the definition of $B_{i_k}$ that
$$\forall\x{A}\in\CX{A}: \qquad \lFac{i_k} \quad\le\quad \frac{\phantom{\sum_{\x{i_k}}} \sum_{\x{i_{k-1}}} \dots \sum_{\x{i_1}} \Psi}{\sum_{\x{i_k}} \sum_{\x{i_{k-1}}} \dots \sum_{\x{i_1}} \Psi} \quad\le\quad \uFac{i_k}$$
for all $k=1,2,\dots,n$.
Taking the product of these $n$ inequalities yields
$$\prod_{k=1}^n \lFac{i_k} \le \nor \Psi \le \prod_{k=1}^n \uFac{i_k}$$
pointwise, and therefore $\nor \Psi \in \sbbox{\prod_{k=1}^n B_{i_k}}$.
\end{proof}

The following corollary is somewhat elaborate to state, but readily follows from
the previous lemma after attaching a factor $I$ that depends on all nodes in $A$
and one additional newly introduced node $i$:

\begin{coro}\label{coro:second_basic_lemma}
Let $\graph = (\vars,\factors,\uedges)$ be a factor graph. Let $i \in \vars$ with exactly
one neighbor in $\factors$, say $\nbv{i} = \{I\}$. Then $\Prob(\x{i}) \in B_i$ where
\begin{equation}\label{eq:bound_on_i_coro_second_basic_lemma}
B_i = \snbbox{\sum_{\x{\nbfe{I}{i}}} \fac{I} \sbbox{\prod_{k \in \nbfe{I}{i}} B_k^{\setm I}}}
\end{equation}
and 
$$B_k^{\setm I} = \snbbox{\sum_{\x{\vars \setm \{i,k\}}} \Psi_{\factors\setm I} \fmeas{\nbfe{I}{\{i,k\}}}^*}$$
with
\begin{equation}
\Psi_{\factors\setm I} := \prod_{J \in \factors\setm I} \psi_J.\qedtag{12pt}
\end{equation}
\end{coro}

We now proceed with a sketch of the proof of Theorem \ref{theo:saw}, which
was inspired by \citep{Ihler07}.\bigskip

\begin{proofsketch}
The proof proceeds using structural induction, recursively transforming 
the original factor graph $\graph$ into the SAW tree $T^{SAW}_{\graph}(i)$,
refining the bound at each step, until it becomes equivalent to the result
of the message propagation algorithm on the SAW tree described above
in  equations \eref{eq:boxprop_saw_cycle_induced_leaf_node}--\eref{eq:boxprop_saw_final_result}.

Let $\graph := (\vars, \factors, \uedges)$ be a factor graph.
Let $i \in \vars$ and let $T^{SAW}_{\graph}(i)$ be the
SAW tree with root $i$. Let $\{I_1,\dots,I_n\} = \nbv{i}$. 

Suppose that $n > 1$. Consider the equivalent factor graph $\graph'$ that is
obtained by creating $n$ copies $i_n$ of the variable node $i$, where each copy
$i_j$ is only connected with the factor node $I_j$ (for $j=1,\dots,n$); in
addition, all copies are connected with the original variable $i$ using the
delta function $\fac{\delta} := \delta(\x{i},\x{i_1},\dots,\x{i_n})$. This step
is illustrated in Figure \ref{fig:saw_proof_var}(a)--(b). Applying Corollary 
\ref{coro:second_basic_lemma} to $\graph'$ yields the following bound which
follows from \eref{eq:bound_on_i_coro_second_basic_lemma} because of the properties of the delta function:
\begin{equation}\label{eq:boundprop_product}
\Prob(\x{i}) \in \snbbox{\prod_{j=1}^n B_{i_j}^{\setm \delta}}
\end{equation}
where
$$B_{i_j}^{\setm \delta} := \snbbox{\sum_{\x{\setm i_j}} \left(\prod_{J\in\factors} \psi_J\right) \fmeas{\{i_1,\dots,i_{j-1},i_{j+1},\dots,i_n\}}^*}\qquad j=1,\dots,n.$$
In the expression on the right-hand side, the factor $\psi_{I_k}$ implicitly depends on $i_k$ instead of $i$ (for all $k=1,\dots,n$).
This bound is represented graphically in Figure \ref{fig:saw_proof_var}(c)--(d)
where the gray variable nodes correspond to simplices of single-variable factors,
\ie they are meant to be multiplied with unknown single-variable factors.
Note that \eref{eq:boundprop_product} corresponds precisely with \eref{eq:boxprop_saw_final_result}.

\begin{figure}[p]
\begin{tikzpicture}
  \tikzstyle{var}=[circle,draw=black,fill=white,semithick,minimum size=15pt]
  \tikzstyle{varg}=[circle,draw=black,fill=gray!20,semithick,minimum size=15pt]
  \tikzstyle{fac}=[rectangle,draw=black,fill=white,semithick,minimum size=12pt]
  \begin{scope}[xshift=0cm,yshift=-0.6cm]
  \small
  \node (a) at (-2.0,1.4) {(a)};
  \begin{scope}[xshift=0.5cm]
  \tiny
  \node (i)  [var] at (0,-0.2){$i$};
  \node (I1) [fac] at (-1,-1) {$I_1$};
  \node (I2) [fac] at (0,-1)  {$I_2$};
  \node (I3) [fac] at (1,-1)  {$I_3$};
  \node (d1)       at (-1,-1.9) {$\vdots$};
  \node (d2)       at (0,-1.9)  {$\vdots$};
  \node (d3)       at (1,-1.9)  {$\vdots$};
  \draw (i) -- (I1);
  \draw (i) -- (I2);
  \draw (i) -- (I3);
  \draw (I1) -- (-1.2,-1.6);
  \draw (I1) -- (-0.8,-1.6);
  \draw (I2) -- (0,-1.6);
  \draw (I3) -- (1.2,-1.6);
  \draw (I3) -- (0.8,-1.6);
  \draw[dotted] (-2,-1.6) rectangle (2,-2.4);
  \end{scope}
  \end{scope}
  \begin{scope}[xshift=6.5cm,yshift=-0.6cm]
  \small
  \node (b) at (-2.0,1.4) {(b)};
  \begin{scope}[xshift=0.5cm]
  \tiny
  \node (i)  [var] at (0,1.4) {$i$};
  \node (delta) [fac] at (0,0.6) {$\delta$};
  \node (i1) [var] at (-1,-0.2)   {$i_1$};
  \node (i2) [var] at (0,-0.2)   {$i_2$};
  \node (i3) [var] at (1,-0.2)   {$i_3$};
  \node (I1) [fac] at (-1,-1) {$I_1$};
  \node (I2) [fac] at (0,-1)  {$I_2$};
  \node (I3) [fac] at (1,-1)  {$I_3$};
  \node (d1)       at (-1,-1.9) {$\vdots$};
  \node (d2)       at (0,-1.9)  {$\vdots$};
  \node (d3)       at (1,-1.9)  {$\vdots$};
  \draw (i) -- (delta);
  \draw (delta) -- (i1);
  \draw (delta) -- (i2);
  \draw (delta) -- (i3);
  \draw (i1) -- (I1);
  \draw (i2) -- (I2);
  \draw (i3) -- (I3);
  \draw (I1) -- (-1.2,-1.6);
  \draw (I1) -- (-0.8,-1.6);
  \draw (I2) -- (0,-1.6);
  \draw (I3) -- (1.2,-1.6);
  \draw (I3) -- (0.8,-1.6);
  \draw[dotted] (-2,-1.6) rectangle (2,-2.4);
  \end{scope}
  \end{scope}
  \begin{scope}[xshift=0cm,yshift=-3.8cm]
  \small
  \node (c) at (-2.0,-0.4) {(c)};
  \tiny
  \begin{scope}[xshift=0cm]
    \colorlet{c1}{gray!0}
    \colorlet{c2}{gray!20}
    \colorlet{c3}{gray!20}
    \node (ia1) [var,fill=c1] at (-1,-2) {$i_{11}$};
    \node (ib1) [var,fill=c2] at (0,-2) {$i_{12}$};
    \node (ic1) [var,fill=c3] at (1,-2) {$i_{13}$};
    \node (I1a) [fac] at (-1,-2.8) {$I_1$};
    \node (I2a) [fac] at (0,-2.8)  {$I_2$};
    \node (I3a) [fac] at (1,-2.8)  {$I_3$};
    \node (d1a)       at (-1,-3.7) {$\vdots$};
    \node (d2a)       at (0,-3.7)  {$\vdots$};
    \node (d3a)       at (1,-3.7)  {$\vdots$};
    \draw (ia1) -- (I1a);
    \draw (ib1) -- (I2a);
    \draw (ic1) -- (I3a);
    \draw (I1a) -- (-1.2,-3.4);
    \draw (I1a) -- (-0.8,-3.4);
    \draw (I2a) -- (0,-3.4);
    \draw (I3a) -- (1.2,-3.4);
    \draw (I3a) -- (0.8,-3.4);
    \draw[dotted] (-1.7,-3.4) rectangle (1.7,-4.2);
  \end{scope}
  \begin{scope}[xshift=4cm]
    \colorlet{c1}{gray!20}
    \colorlet{c2}{gray!0}
    \colorlet{c3}{gray!20}
    \node (ia2) [var,fill=c1] at (-1,-2) {$i_{21}$};
    \node (ib2) [var,fill=c2] at (0,-2) {$i_{22}$};
    \node (ic2) [var,fill=c3] at (1,-2) {$i_{23}$};
    \node (I1a) [fac] at (-1,-2.8) {$I_1$};
    \node (I2a) [fac] at (0,-2.8)  {$I_2$};
    \node (I3a) [fac] at (1,-2.8)  {$I_3$};
    \node (d1a)       at (-1,-3.7) {$\vdots$};
    \node (d2a)       at (0,-3.7)  {$\vdots$};
    \node (d3a)       at (1,-3.7)  {$\vdots$};
    \draw (ia2) -- (I1a);
    \draw (ib2) -- (I2a);
    \draw (ic2) -- (I3a);
    \draw (I1a) -- (-1.2,-3.4);
    \draw (I1a) -- (-0.8,-3.4);
    \draw (I2a) -- (0,-3.4);
    \draw (I3a) -- (1.2,-3.4);
    \draw (I3a) -- (0.8,-3.4);
    \draw[dotted] (-1.7,-3.4) rectangle (1.7,-4.2);
  \end{scope}
  \begin{scope}[xshift=8cm]
    \colorlet{c1}{gray!20}
    \colorlet{c2}{gray!20}
    \colorlet{c3}{gray!0}
    \node (ia3) [var,fill=c1] at (-1,-2) {$i_{31}$};
    \node (ib3) [var,fill=c2] at (0,-2) {$i_{32}$};
    \node (ic3) [var,fill=c3] at (1,-2) {$i_{33}$};
    \node (I1a) [fac] at (-1,-2.8) {$I_1$};
    \node (I2a) [fac] at (0,-2.8)  {$I_2$};
    \node (I3a) [fac] at (1,-2.8)  {$I_3$};
    \node (d1a)       at (-1,-3.7) {$\vdots$};
    \node (d2a)       at (0,-3.7)  {$\vdots$};
    \node (d3a)       at (1,-3.7)  {$\vdots$};
    \draw (ia3) -- (I1a);
    \draw (ib3) -- (I2a);
    \draw (ic3) -- (I3a);
    \draw (I1a) -- (-1.2,-3.4);
    \draw (I1a) -- (-0.8,-3.4);
    \draw (I2a) -- (0,-3.4);
    \draw (I3a) -- (1.2,-3.4);
    \draw (I3a) -- (0.8,-3.4);
    \draw[dotted] (-1.7,-3.4) rectangle (1.7,-4.2);
  \end{scope}
  \node (delta) [fac] at (4,-1.2) {$\delta$};
  \node (i)     [var] at (4,-0.4)  {$i$};
  \draw (i.south) -- (delta.north);
  \draw (delta.west) -|  (ia1.north);
  \draw (delta.south) -- (ib2.north);
  \draw (delta.east) -|  (ic3.north);
  \end{scope}
  \begin{scope}[xshift=0cm,yshift=-8cm]
  \small
  \node (d) at (-2.0,-1.2) {(d)};
  \tiny
  \begin{scope}[xshift=0cm]
  \node (Ia1) [fac] at (0,-2) {$I_1$};
  \node             at (0,-2.9) {$\vdots$};
  \node (Ib1) [fac] at (-0.5,-4) {$I_2$};
  \node (Ic1) [fac] at (0.5,-4) {$I_3$};
  \node (ib1) [varg] at (-0.5,-4.8) {$i_{12}$};
  \node (ic1) [varg] at (0.5,-4.8) {$i_{13}$};
  \draw (Ib1) -- (ib1);
  \draw (Ic1) -- (ic1);
  \draw (Ia1) -- (-.2,-2.6);
  \draw (Ia1) -- (0.2,-2.6);
  \draw (Ib1) -- (-0.5,-3.4);
  \draw (Ic1) -- (0.3,-3.4);
  \draw (Ic1) -- (0.7,-3.4);
  \draw[dotted] (-1.7,-2.6) rectangle (1.7,-3.4);
  \end{scope}
  \begin{scope}[xshift=4cm]
  \node (Ib2) [fac] at (0,-2) {$I_2$};
  \node             at (0,-2.9) {$\vdots$};
  \node (Ia2) [fac] at (-0.5,-4) {$I_1$};
  \node (Ic2) [fac] at (0.5,-4) {$I_3$};
  \node (ia2) [varg] at (-0.5,-4.8) {$i_{21}$};
  \node (ic2) [varg] at (0.5,-4.8) {$i_{23}$};
  \draw (Ia2) -- (ia2);
  \draw (Ic2) -- (ic2);
  \draw (Ib2) -- (0,-2.6);
  \draw (Ia2) -- (-0.7,-3.4);
  \draw (Ia2) -- (-0.3,-3.4);
  \draw (Ic2) -- (0.3,-3.4);
  \draw (Ic2) -- (0.7,-3.4);
  \draw[dotted] (-1.7,-2.6) rectangle (1.7,-3.4);
  \end{scope}
  \begin{scope}[xshift=8cm]
  \node (Ic3) [fac] at (0,-2) {$I_3$};
  \node             at (0,-2.9) {$\vdots$};
  \node (Ia3) [fac] at (-0.5,-4) {$I_1$};
  \node (Ib3) [fac] at (0.5,-4) {$I_2$};
  \node (ia3) [varg] at (-0.5,-4.8) {$i_{31}$};
  \node (ib3) [varg] at (0.5,-4.8) {$i_{32}$};
  \draw (Ia3) -- (ia3);
  \draw (Ib3) -- (ib3);
  \draw (Ic3) -- (-0.2,-2.6);
  \draw (Ic3) -- (0.2,-2.6);
  \draw (Ia3) -- (-0.7,-3.4);
  \draw (Ia3) -- (-0.3,-3.4);
  \draw (Ib3) -- (0.5,-3.4);
  \draw[dotted] (-1.7,-2.6) rectangle (1.7,-3.4);
  \end{scope}
  \node (i)     [var] at (4,-1.2)  {$i$};
  \draw (i) -- (Ia1.north);
  \draw (i) -- (Ib2.north);
  \draw (i) -- (Ic3.north);
  \end{scope}
\end{tikzpicture}
\caption{\label{fig:saw_proof_var}One step in the proof of Theorem \ref{theo:saw}: propagating 
bounds towards variable $i$ in case it has more than one neighboring factor nodes 
$I_1, \dots, I_n$ (here, $n=3$). Gray nodes represent added (unknown) single-variable factors.
(a) Factor graph $\graph$. (b) Equivalent factor graph $\graph'$.
(c) Result of replicating $\graph$ $n$ times, where in each 
copy $\graph_k$ of $\graph$, $i$ is replaced by exactly $n$ copies $i_{kj}$ of $i$
for $j=1,\dots,n$, where $i_{kj}$ is connected only with the factor $I_j$ in $\graph_k$.
Then, the original variable $i$ is connected using a delta factor with $n$ of its copies
$i_{jj}$ for $j=1,\dots,n$.
(d) Simplification of (c) obtained by identifying $i$ with its $n$ copies $i_{jj}$ for 
$j=1,\dots,n$ and changing the layout slightly.}
\end{figure}
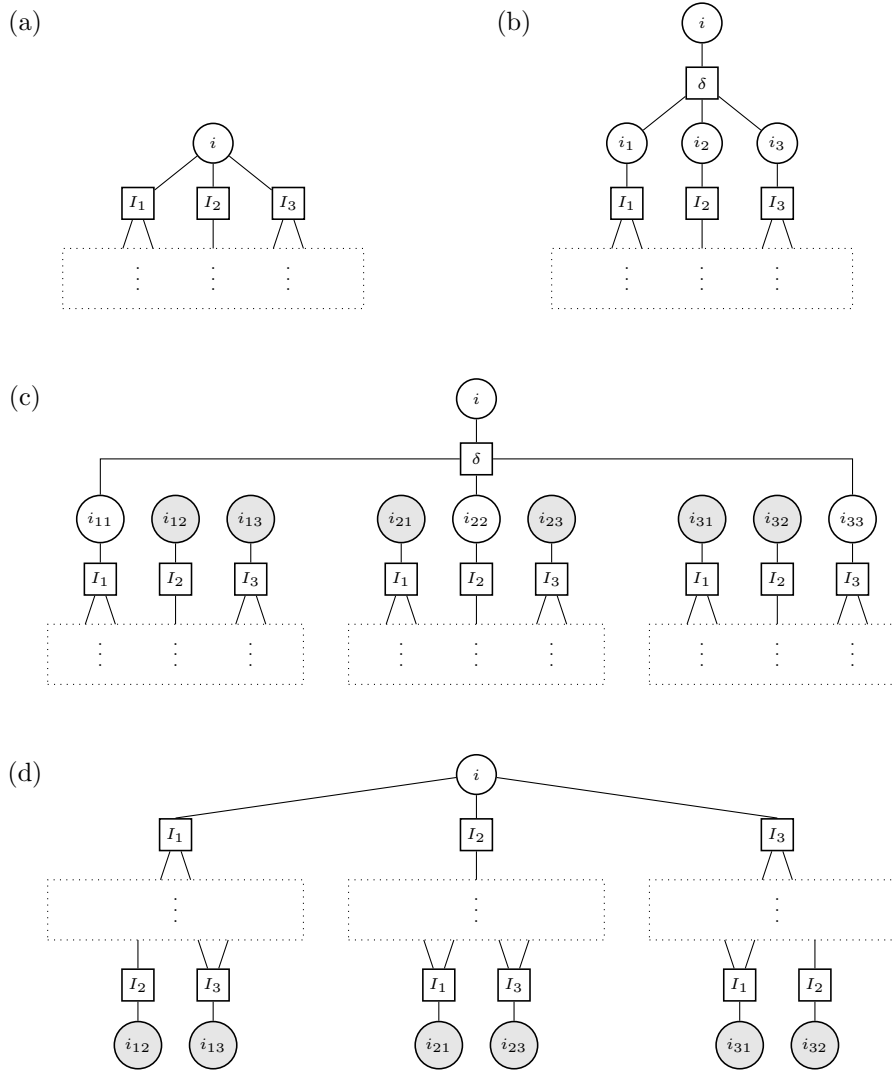

\begin{figure}[p]
\begin{tikzpicture}
  \tikzstyle{var}=[circle,draw=black,fill=white,semithick,minimum size=15pt]
  \tikzstyle{varg}=[circle,draw=black,fill=gray!20,semithick,minimum size=15pt]
  \tikzstyle{fac}=[rectangle,draw=black,fill=white,semithick,minimum size=12pt]
  \tikzstyle{facg}=[rectangle,draw=black,fill=gray!20,semithick,minimum size=12pt]
  \begin{scope}[xshift=0cm,yshift=0cm]
  \small
  \node (a) at (-2.0,1) {(a)};
  \begin{scope}[xshift=4cm]
  \tiny
  \node (i)  [var] at (0,1)   {$i$};
  \node (I)  [fac] at (0,0)   {$I$};
  \node (j1) [var] at (-2,-1) {$j_1$};
  \node (j2) [var] at (0,-1)  {$j_2$};
  \node (j3) [var] at (2,-1)  {$j_3$};
  \node            at (-2,-2) {$\vdots$};
  \node            at (0,-2)  {$\vdots$};
  \node            at (2,-2)  {$\vdots$};
  \draw (i) -- (I);
  \draw (I) -- (j1);
  \draw (I) -- (j2);
  \draw (I) -- (j3);
  \draw (j1) -- (-2.2,-1.7);
  \draw (j1) -- (-1.8,-1.7);
  \draw (j2) -- (0,-1.7);
  \draw (j3) -- (2.2,-1.7);
  \draw (j3) -- (1.8,-1.7);
  \draw[dotted] (-3,-1.7) rectangle (3,-2.6);
  \end{scope}
  \end{scope}
  \begin{scope}[xshift=0cm,yshift=-3.7cm]
  \small
  \node (b) at (-2.0,0) {(b)};
  \tiny
  \begin{scope}[xshift=0cm]
    \colorlet{c1}{gray!0}
    \colorlet{c2}{gray!20}
    \colorlet{c3}{gray!20}
    \node (ja1) [var,fill=c1] at (-1,-2) {$j_{11}$};
    \node (jb1) [var,fill=c2] at (0,-2) {$j_{12}$};
    \node (jc1) [var,fill=c3] at (1,-2) {$j_{13}$};
    \node             at (-1,-3) {$\vdots$};
    \node             at (0,-3)  {$\vdots$};
    \node             at (1,-3)  {$\vdots$};
    \draw (ja1) -- (-1.2,-2.7);
    \draw (ja1) -- (-0.8,-2.7);
    \draw (jb1) -- (0,-2.7);
    \draw (jc1) -- (1.2,-2.7);
    \draw (jc1) -- (0.8,-2.7);
    \draw[dotted] (-1.7,-2.7) rectangle (1.7,-3.6);
  \end{scope}
  \begin{scope}[xshift=4cm]
    \colorlet{c1}{gray!20}
    \colorlet{c2}{gray!0}
    \colorlet{c3}{gray!20}
    \node (ja2) [var,fill=c1] at (-1,-2) {$j_{21}$};
    \node (jb2) [var,fill=c2] at (0,-2) {$j_{22}$};
    \node (jc2) [var,fill=c3] at (1,-2) {$j_{23}$};
    \node             at (-1,-3) {$\vdots$};
    \node             at (0,-3)  {$\vdots$};
    \node             at (1,-3)  {$\vdots$};
    \draw (ja2) -- (-1.2,-2.7);
    \draw (ja2) -- (-0.8,-2.7);
    \draw (jb2) -- (0,-2.7);
    \draw (jc2) -- (1.2,-2.7);
    \draw (jc2) -- (0.8,-2.7);
    \draw[dotted] (-1.7,-2.7) rectangle (1.7,-3.6);
  \end{scope}
  \begin{scope}[xshift=8cm]
    \colorlet{c1}{gray!20}
    \colorlet{c2}{gray!20}
    \colorlet{c3}{gray!0}
    \node (ja3) [var,fill=c1] at (-1,-2) {$j_{31}$};
    \node (jb3) [var,fill=c2] at (0,-2) {$j_{32}$};
    \node (jc3) [var,fill=c3] at (1,-2) {$j_{33}$};
    \node             at (-1,-3) {$\vdots$};
    \node             at (0,-3)  {$\vdots$};
    \node             at (1,-3)  {$\vdots$};
    \draw (ja3) -- (-1.2,-2.7);
    \draw (ja3) -- (-0.8,-2.7);
    \draw (jb3) -- (0,-2.7);
    \draw (jc3) -- (1.2,-2.7);
    \draw (jc3) -- (0.8,-2.7);
    \draw[dotted] (-1.7,-2.7) rectangle (1.7,-3.6);
  \end{scope}
  \node (I) [fac] at (4,-1) {$I$};
  \node (i) [var] at (4,0)  {$i$};
  \draw (i.south) -- (I.north);
  \draw (I.west) -|  (ja1.north);
  \draw (I.south) -- (jb2.north);
  \draw (I.east) -|  (jc3.north);
  \end{scope}
  \begin{scope}[xshift=0cm,yshift=-8.4cm]
  \small
  \node (c) at (-2.0,0) {(c)};
  \tiny
  \begin{scope}[xshift=0cm]
  \node (ja1) [var] at (0,-2) {$j_{11}$};
  \node             at (0,-3) {$\vdots$};
  \node (jb1) [var] at (-0.5,-4.3) {$j_{12}$};
  \node (jc1) [var] at (0.5,-4.3) {$j_{13}$};
  \node (Ib1) [facg] at (-0.5,-5.3) {$I$};
  \node (Ic1) [facg] at (0.5,-5.3) {$I$};
  \draw (jb1) -- (Ib1);
  \draw (jc1) -- (Ic1);
  \draw (ja1) -- (-.2,-2.7);
  \draw (ja1) -- (0.2,-2.7);
  \draw (jb1) -- (-0.5,-3.6);
  \draw (jc1) -- (0.3,-3.6);
  \draw (jc1) -- (0.7,-3.6);
  \draw[dotted] (-1.7,-2.7) rectangle (1.7,-3.6);
  \end{scope}
  \begin{scope}[xshift=4cm]
  \node (jb2) [var]  at (0,-2) {$j_{22}$};
  \node              at (0,-3) {$\vdots$};
  \node (ja2) [var] at (-0.5,-4.3) {$j_{21}$};
  \node (jc2) [var] at (0.5,-4.3) {$j_{23}$};
  \node (Ia2) [facg] at (-0.5,-5.3) {$I$};
  \node (Ic2) [facg] at (0.5,-5.3) {$I$};
  \draw (ja2) -- (Ia2);
  \draw (jc2) -- (Ic2);
  \draw (jb2) -- (0,-2.7);
  \draw (ja2) -- (-0.7,-3.6);
  \draw (ja2) -- (-0.3,-3.6);
  \draw (jc2) -- (0.3,-3.6);
  \draw (jc2) -- (0.7,-3.6);
  \draw[dotted] (-1.7,-2.7) rectangle (1.7,-3.6);
  \end{scope}
  \begin{scope}[xshift=8cm]
  \node (jc3) [var] at (0,-2) {$j_{33}$};
  \node             at (0,-3) {$\vdots$};
  \node (ja3) [var] at (-0.5,-4.3) {$j_{31}$};
  \node (jb3) [var] at (0.5,-4.3) {$j_{32}$};
  \node (Ia3) [facg] at (-0.5,-5.3) {$I$};
  \node (Ib3) [facg] at (0.5,-5.3) {$I$};
  \draw (ja3) -- (Ia3);
  \draw (jb3) -- (Ib3);
  \draw (jc3) -- (-0.2,-2.7);
  \draw (jc3) -- (0.2,-2.7);
  \draw (ja3) -- (-0.7,-3.6);
  \draw (ja3) -- (-0.3,-3.6);
  \draw (jb3) -- (0.5,-3.6);
  \draw[dotted] (-1.7,-2.7) rectangle (1.7,-3.6);
  \end{scope}
  \node (i)     [var] at (4,0)  {$i$};
  \node (I)     [fac] at (4,-1) {$I$};
  \draw (i) -- (I);
  \draw (I) -- (ja1.north);
  \draw (I) -- (jb2.north);
  \draw (I) -- (jc3.north);
  \end{scope}
\end{tikzpicture}
\caption{\label{fig:saw_proof_fac}Another step in the proof of Theorem \ref{theo:saw}: propagating
bounds towards variable $i$ in case it has exactly one neighboring factor node $I$ which has $m+1$
neighboring variables $\{i,j_1,\dots,j_m\}$. (a) Factor graph $\graph$. (b) Result of replicating 
$\graph\setm\{i,I\}$ $m$ times and connecting the factor $I$ with $i$ and with copy $j_{kk}$ of $j_k$
for $k=1,\dots,m$. (c) Equivalent to (b) but with a slightly changed layout. The gray copies of $I$
represent (unknown) single-variable factors (on their neighboring variable).}
\end{figure}
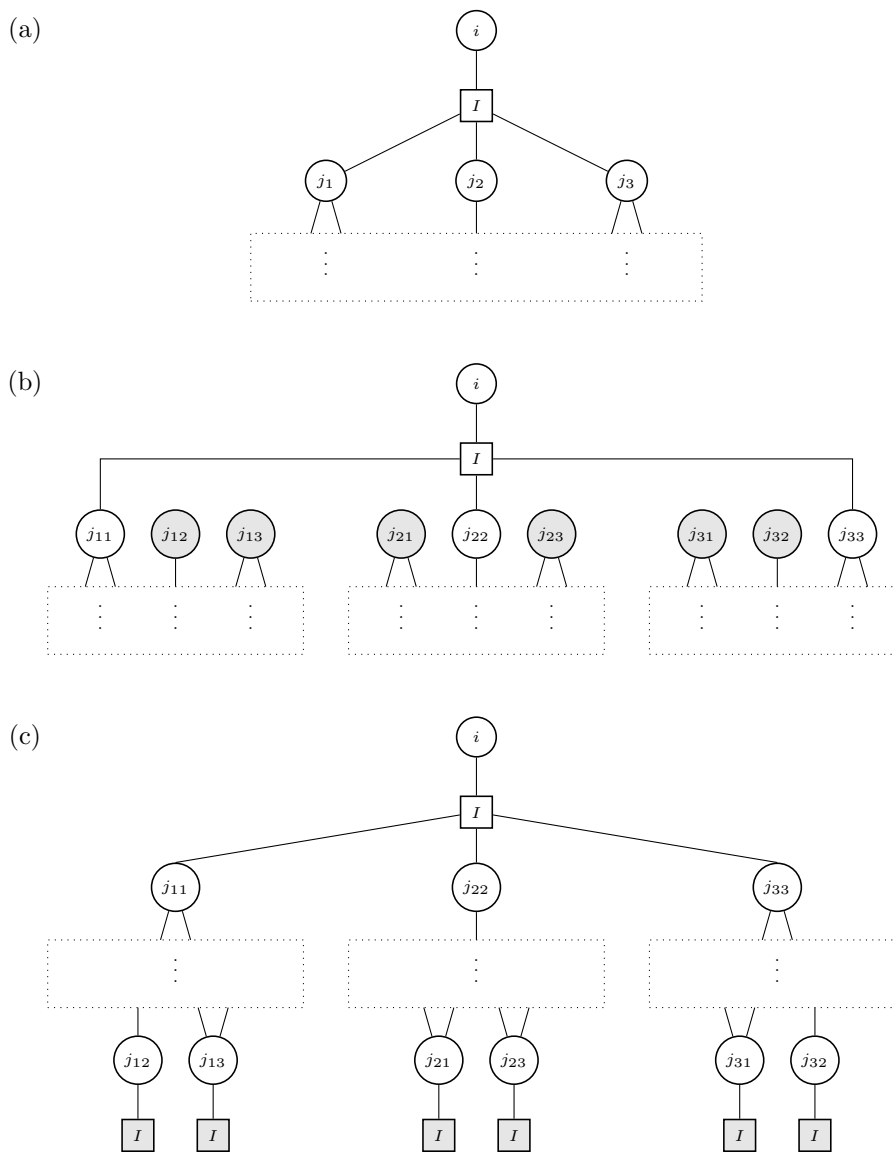

\clearpage

The case that $n = 1$ is simpler because there is no need to introduce the delta function.
It is illustrated in Figure \ref{fig:saw_proof_fac}. Let $\{I\} = \nbv{i}$ and let $\{j_1,\dots,j_m\} = \nbfe{I}{i}$.
Applying Corollary \ref{coro:second_basic_lemma} yields the following bound:
\begin{equation}\label{eq:boundprop_prodsum}
\Prob(\x{i}) \in \snbbox{\sum_{\x{\nbfe{I}{i}}} \fac{I} \sbbox{\prod_{k=1}^m B_{j_k}^{\setm I}}}
\end{equation}
where
$$B_{j_k}^{\setm I} := \snbbox{\sum_{\x{\setm \{i,j_k\}}} \left(\prod_{J \in \factors\setm I} \psi_J\right) \fmeas{\{j_1,\dots,j_{k-1},j_{k+1},\dots,j_m\}}^*}\qquad k=1,\dots,m.$$
This bound is represented graphically in Figure \ref{fig:saw_proof_fac}(b)--(c)
where the gray nodes correspond with simplices of single-variable factors.
Note that \eref{eq:boundprop_prodsum} corresponds precisely with \eref{eq:boxprop_saw_factor_to_variable}.

Recursively iterating the factor graph operations in Figures \ref{fig:saw_proof_fac} and 
\ref{fig:saw_proof_var}, the connected component that remains in the end is precisely the 
SAW tree $T^{SAW}_{\graph}(i)$; the bounds derived along the way correspond precisely with 
the message passing algorithm on the SAW tree described above.
\end{proofsketch}

Again, the self-avoiding walk tree with root $i$ is a subtree of the computation
tree for $i$. This means that the bounds on the
exact marginals given by Theorem \ref{theo:saw} are bounds on the approximate
Belief Propagation marginals (beliefs) as well.

\begin{coro}\label{coro:saw_bp}
In the situation described in Theorem \ref{theo:saw}, the final bounding
box $\bebox{i}$ also bounds the (approximate) Belief Propagation marginal
of the root node $i$, \ie $\Prob_{BP}(\x{i}) \in \bebox{i}$.\qed
\end{coro}

\section{Related work}\label{sec5:related_work}

There exist many other bounds on single-variable marginals. Also,
bounds on the partition sum can be used to obtain bounds on
single-variable marginals. For all bounds known to the authors,
we will discuss how they compare with our bounds.
In the following, we will denote exact marginals as
$p_i(\x{i}) := \Prob(\x{i})$ and BP marginals (beliefs) as $b_i(\x{i}) := \Prob_{BP}(\x{i})$.
%

\subsection{The Dobrushin-Tatikonda bound}
\citet{Tatikonda03} derived a bound on the error of BP marginals using
mathematical tools from Gibbs measure theory \citep{Georgii88}, in particular
using a result known as Dobrushin's theorem. The bounds on the error of the
BP marginals can be easily translated into bounds on the exact marginals:
$$\abs{b_i(\x{i}) - p_i(\x{i})} \le \epsilon \implies p_i(\x{i}) \in [b_i(\x{i})-\epsilon,b_i(\x{i})+\epsilon]$$
for all $i\in\vars$ and $\x{i}\in\CX{i}$.

The Dobrushin-Tatikonda bound depends on the \emph{girth} of the graph (the
number of edges in the shortest cycle, or infinity if there is no cycle) and
the properties of Dobrushin's interdependence matrix, which is a $N \times N$
matrix $C$. The entry $C_{ij}$ is only nonzero if $i \in \del{j}$ and in that
case, the computational cost of computing its value is exponential in the size
of the Markov blanket. Thus the computational complexity of the
Dobrushin-Tatikonda bound is $\C{O} (\max_{i\in\vars} \nel{\CX{\del{i}}})$,
plus the cost of running BP.

\subsection{The Dobrushin-Taga-Mase bound}
Inspired by the work of \citet{TatikondaJordan02}, \citet{TagaMase06b} derived another 
bound on the error of BP marginals, also based on Dobrushin's theorem. 
This bound also depends on the properties of Dobrushin's interdependence matrix 
and has similar computational cost. Whereas the Dobrushin-Tatikonda bound gives
one bound for all variables, the Dobrushin-Taga-Mase bound gives a different
bound for each variable.

\subsection{Bound Propagation}
\citet{LeisinkKappen03} proposed a method called ``Bound Propagation'' 
which can be used to obtain bounds on exact marginals. The idea underlying 
this method is very similar to the one employed in this work, with one
crucial difference. Whereas we use a cavity approach, using as basis equation
\begin{equation*}
\Prob(\x{i}) \propto \sum_{\xd{i}} \left(\prod_{I \in \nbv{i}} \fac{I}\right) \Pmx{i},
\qquad \Pmx{i} \propto \sum_{\x{\vars \setm \Del{i}}} \prod_{I\in\factors\setm \nbf{i}} \fac{I}
\end{equation*}
and bound the quantity $\Prob(\x{i})$ by optimizing over $\Pmx{i}$, the 
basis equation employed by Bound Propagation is
\begin{equation*}
\Prob(\x{i}) = \sum_{\xd{i}} \Prob(\x{i} \given \xd{i}) \Prob(\xd{i})
\end{equation*}
and the optimization is over $\Prob(\xd{i})$. Unlike in our case,
the computational complexity is exponential in the size of the
Markov blanket, because of the required calculation of the conditional
distribution $\Prob(\x{i} \given \xd{i})$. On the other hand, the advantage
of this approach is that a bound on $\Prob(\x{j})$ for $j \in \del{i}$ is also
a bound on $\Prob(\xd{i})$, which in turn gives rise to a bound on $\Prob(\x{i})$.
In this way, bounds can propagate through the graphical
model, eventually yielding a new (tighter) bound on $\Prob(\xd{i})$. 
Although the iteration can result in rather tight bounds, the main 
disadvantage of Bound Propagation is its computational cost: it is 
exponential in the Markov blanket and often many iterations are 
needed for the bounds to become tight. Indeed, for a 
simple tree of $N = 100$ variables, it can happen that Bound Propagation 
needs several minutes and still obtains very loose bounds (whereas our
bounds give the exact marginal as lower and as upper bound, \ie they
arrive at the optimally tight bound).

\subsection{Upper and lower bounds on the partition sum}
Various upper and lower bounds on the partition sum $Z$ in \eref{eq:prob_dist} exist.
An upper and a lower bound on $Z$ can be combined to obtain bounds
on marginals in the following way. First, note that the exact
marginal of $i$ satisfies
$$\Prob(\x{i}) = \frac{Z_i(\x{i})}{Z},$$
where we defined the partition sum of the \emph{clamped} model as follows:
$$Z_i(\x{i}) := \sum_{\x{\vars\setm \{i\}}} \prod_{I\in\factors} \fac{I}.$$
Thus, we can bound
$$\frac{Z_i^-(\x{i})}{Z^+} \le p_i(\x{i}) \le \frac{Z_i^+(\x{i})}{Z^-}$$
where $Z^- \le Z \le Z^+$ and $Z_i^-(\x{i}) \le Z_i(\x{i}) \le Z_i^+(\x{i})$
for all $\x{i} \in \CX{i}$.

A well-known lower bound of the partition sum is the Mean Field bound. A
tighter lower bound was derived by \citet{LeisinkKappen01}. An upper bound on
the log partition sum was derived by \citet{WainwrightJaakkolaWillsky05}.
Other lower and upper bounds (for the case of binary variables with pairwise
interactions) have been derived by \citet{JaakkolaJordan96}.  

\section{Experiments}\label{sec5:experiments}

\newcommand{\BP}{{\sc BP}\relax}
\newcommand{\BoxPropSubtree}{{\sc BoxProp-SubT}\relax}
\newcommand{\BoxPropSubtreeI}{{\sc Ihler-SubT}\relax}
\newcommand{\BoxPropSAW}{{\sc BoxProp-SAWT}\relax}
\newcommand{\BoxPropSAWI}{{\sc Ihler-SAWT}\relax}
\newcommand{\BoundProp}{{\sc BoundProp}\relax}
\newcommand{\DT}{{\sc DT}\relax}
\newcommand{\DTM}{{\sc DTM}\relax}
\newcommand{\MF}{{\sc MF}\relax}
\newcommand{\LK}{{\sc LK3}\relax}
\newcommand{\JJ}{{\sc JJ}\relax}
\newcommand{\TRW}{{\sc TRW}\relax}

We have done several experiments to compare the quality and computation
time of various bounds empirically. For each variable in the factor graph
under consideration, we calculated the \emph{gap} for each bound 
$\bbox{i}{\lFac{i}}{\uFac{i}} \ni \Prob(\x{i})$, 
which we define as the $\ell_0$-norm $\norms{\uFac{i} - \lFac{i}}{0} = \max_{\x{i}\in\CX{i}} \abs{\uFac{i}(\x{i})-\lFac{i}(\x{i})}$. 

We have used the following bounds in our comparison:
\begin{description}
\item[\DT:] Dobrushin-Tatikonda \citep[Proposition V.6]{Tatikonda03}.
\item[\DTM:] Dobrushin-Taga-Mase \citep[Theorem 1]{TagaMase06b}.
\item[\BoundProp:] Bound Propagation \citep{LeisinkKappen03}, using the implementation of M.\ Leisink, where we chose the maximum cluster size to be $\max_{i\in\vars} \nel{\Del{i}}$.
\item[\BoxPropSubtree:] Theorem \ref{theo:main}, where we used a simple breadth-first algorithm to recursively construct the subtree.
\item[\BoxPropSAW:] Theorem \ref{theo:saw}, where we truncated the SAW tree to at most 5000 nodes.
\item[\BoxPropSAWI:] Ihler's bound \citep{Ihler07}. This bound has only been formulated for pairwise interactions.
\item[\BoxPropSubtreeI:] Ihler's bound \citep{Ihler07} applied on a truncated version of the SAW tree, namely on
the same subtree as used in \BoxPropSubtree. This bound has only been formulated for pairwise interactions.
\end{description}
In addition, we compared with appropriate combinations of the following bounds:
\begin{description}
\item[\MF:] Mean-field lower bound.
\item[\LK:] Third-order lower bound \citep[Eq.\ (10)]{LeisinkKappen01}, where we took for $\mu_i$ the mean field solutions. This bound has been formulated only for the binary, pairwise case.
\item[\JJ:] Refined upper bound \citep[Section 2.2]{JaakkolaJordan96}, with a greedy optimization over the parameters. This bound has been formulated only for the binary, pairwise case.
\item[\TRW:] Our implementation of \citep{WainwrightJaakkolaWillsky05}. This bound has been formulated only for pairwise interactions.
\end{description}
For reference, we calculated the Belief Propagation (\BP) errors by comparing with the exact marginals, using the $\ell_0$ distance as error measure.

\subsection{Grids with binary variables}

We considered a $5\times 5$ Ising grid with binary ($\pm 1$-valued) variables,
i.i.d.\ spin-glass nearest-neighbor interactions $J_{ij} \sim \C{N}(0,\beta^2)$ and i.i.d.\ local fields 
$\theta_i \sim \C{N}(0,\beta^2)$, with probability distribution
\begin{equation*}
\Prob(x) = \frac{1}{Z} \exp\left( \sum_{i\in\vars} \theta_i \x{i} + \frac{1}{2} \sum_{i\in\vars}\sum_{j\in\del{i}} J_{ij} \x{i} \x{j} \right).
\end{equation*}
We took one random instance of the parameters $J$ and $\theta$ (drawn for
$\beta = 1$) and scaled these parameters with the interaction strength
parameter $\beta$, for which we took values in $\{10^{-2},10^{-1},1,10\}$.

\begin{figure}
\centering
\includegraphics[width=0.48\textwidth]{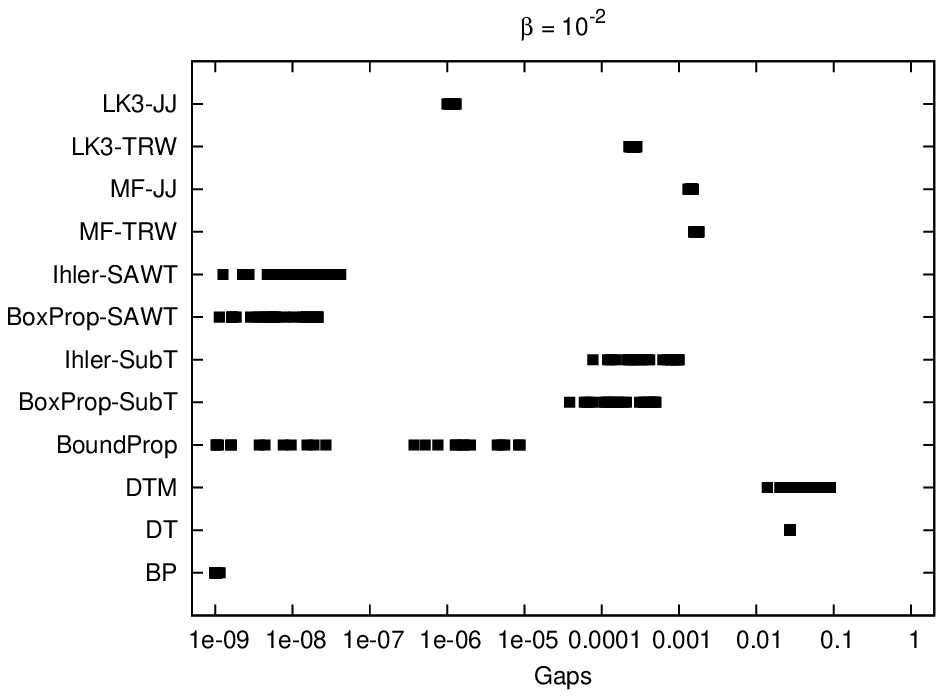}
\includegraphics[width=0.48\textwidth]{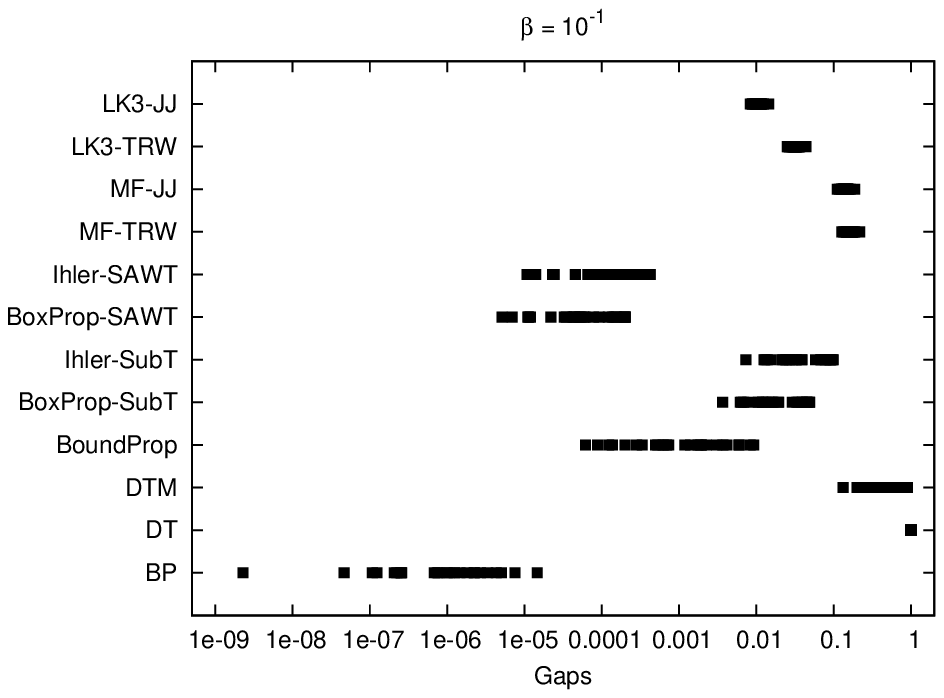}\\[10pt]
\includegraphics[width=0.48\textwidth]{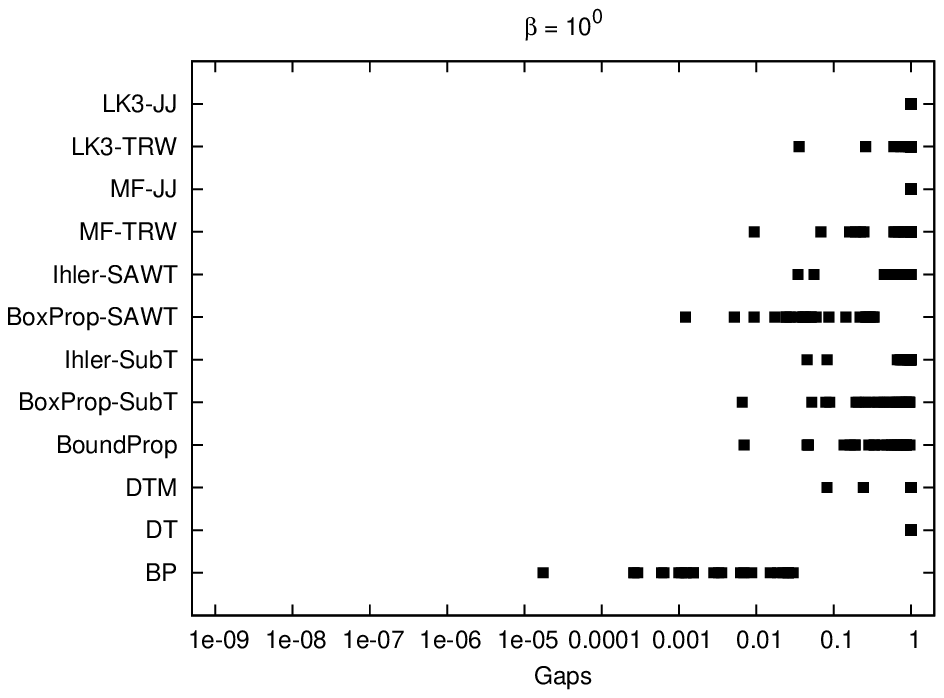}
\includegraphics[width=0.48\textwidth]{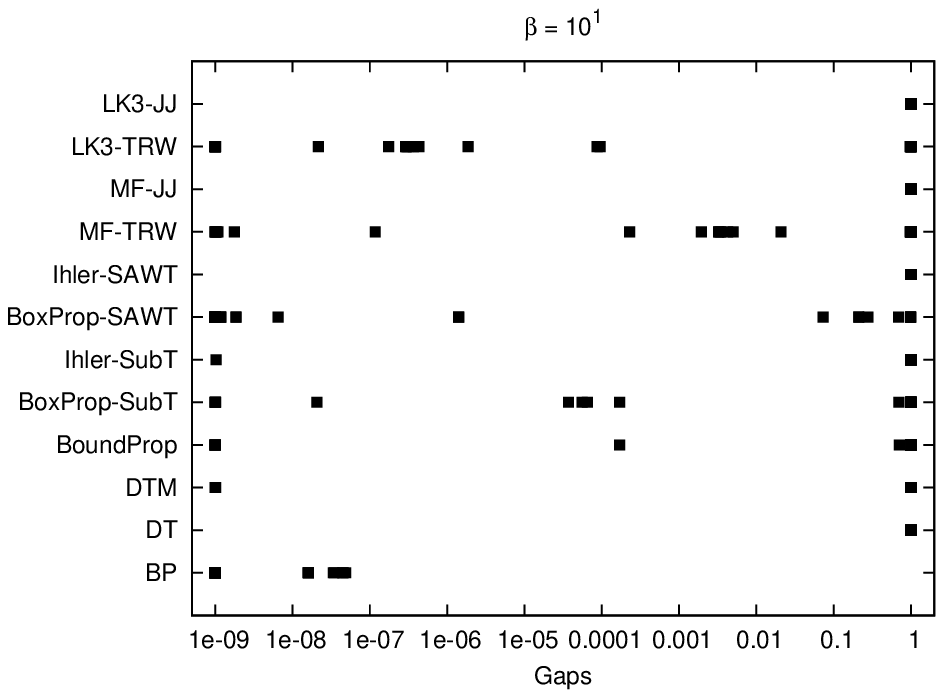}\\[10pt]
\includegraphics[scale=0.7]{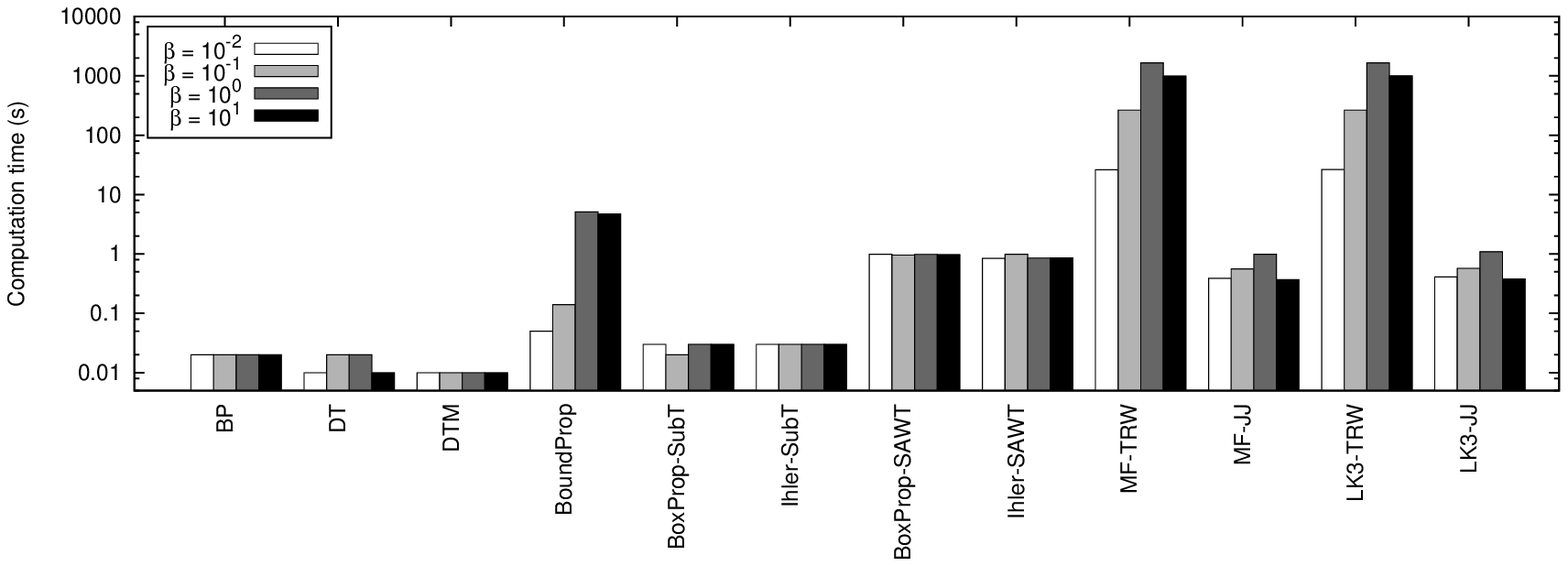}
\caption{Results for grids with binary variables. The first four graphs show, for different values of the interaction strength $\beta$, the gaps of bounds on 
marginals calculated using different methods. Also shown for reference are the errors in the BP 
approximations to the same marginals. The final graph shows the total computation time for each
method.\label{fig:bounds_new}}
\end{figure}

\begin{figure}
\centering
\includegraphics[width=0.3\textwidth]{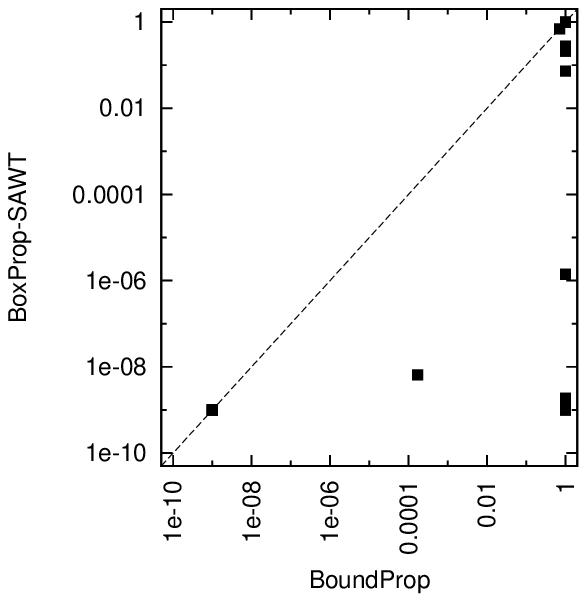}\hfill
\includegraphics[width=0.3\textwidth]{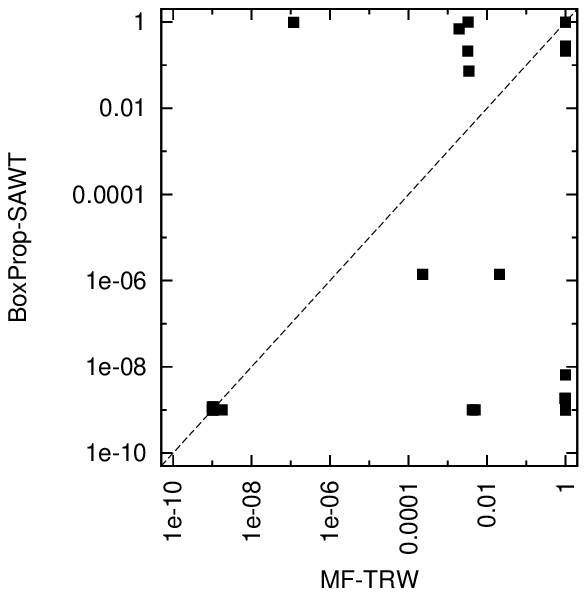}\hfill
\includegraphics[width=0.3\textwidth]{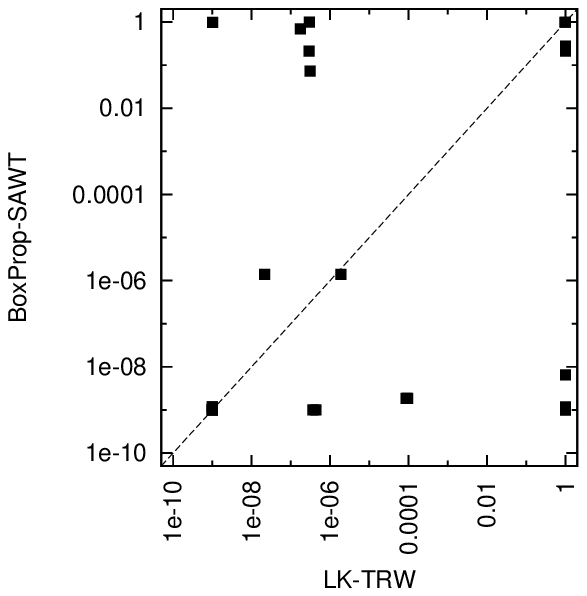}
\caption{Scatter plots comparing some methods in detail for grids with binary variables for strong
interactions ($\beta=10$).\label{fig:bounds_new_scatter}}
\end{figure}

The results are shown in Figure \ref{fig:bounds_new}. 
For very weak interactions ($\beta=10^{-2}$), \BoxPropSAW\ gave the tightest bounds
of all other methods, the only exception being \BoundProp, which gave
a somewhat tighter bound for 5 variables out of 25. For weak and moderate
interactions ($\beta=10^{-1}, 1$), \BoxPropSAW\ gave the tightest bound of all
methods for each variable. For strong interactions ($\beta=10$), the results
were mixed, the best methods being \BoxPropSAW, \BoundProp, \MF-\TRW\ and
\LK-\TRW. Of these, \BoxPropSAW\ was the fastest method, whereas the methods using
\TRW\ were the slowest.\footnote{We had to loosen the convergence criterion for the 
inner loop of \TRW, otherwise it would have taken hours. Since some of the bounds
are significantly tighter than the convergence criterion we used, this may suggest 
that one can loosen the convergence criterion for \TRW\ even more and still obtain
good results using less computation time than the results we present here. 
Unfortunately, it is not clear how this criterion should be chosen in an 
optimal way without actually trying different values and using the best one.} 
For $\beta=10$, we present scatter plots in Figure \ref{fig:bounds_new_scatter} to 
compare the results of some methods in more detail. These plots illustrate that
the tightness of bounds can vary widely over methods and variables.

Among the methods yielding the tightest bounds, the computation time of our bounds
is relatively low in general; only for low interaction strengths \BoundProp\ is 
faster than \BoxPropSAW. Furthermore, the computation time of our bounds does not depend on the interaction strength, 
in contrast with iterative methods such as \BoundProp\ and \MF-\TRW, which need more 
iterations for increasing interaction strength (as the variables become more and more
correlated). Further, as expected, \BoxPropSubtree\ needs less computation time than \BoxPropSAW\ 
but also yields less tight bounds. Another observation is that our bounds
outperform the related versions of Ihler's bounds.

\subsection{Grids with ternary variables}

To evaluate the bounds beyond the special case of binary variables, we have
performed experiments on a $5\times 5$ grid with ternary variables and pairwise
factors between nearest-neighbor variables on the grid. The entries of the 
factors were i.i.d., drawn by taking a random number from a normal distribution
$\normaldist{0}{\beta}$ with mean $0$ and standard deviation $\beta$ and taking
the $\exp$ of that random number.

\begin{figure}
\centering
\includegraphics[width=0.48\textwidth]{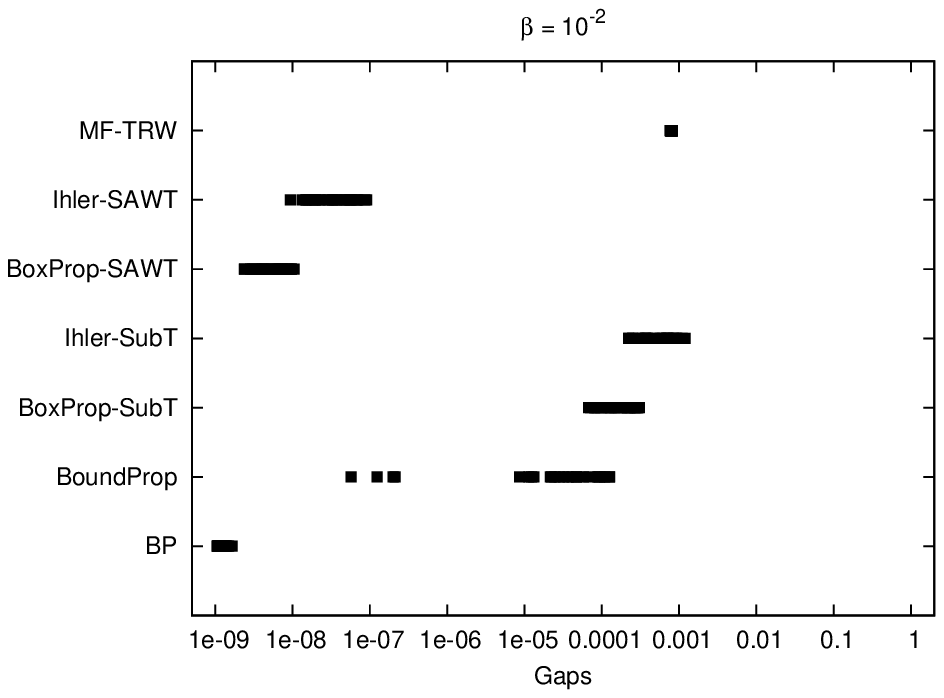}
\includegraphics[width=0.48\textwidth]{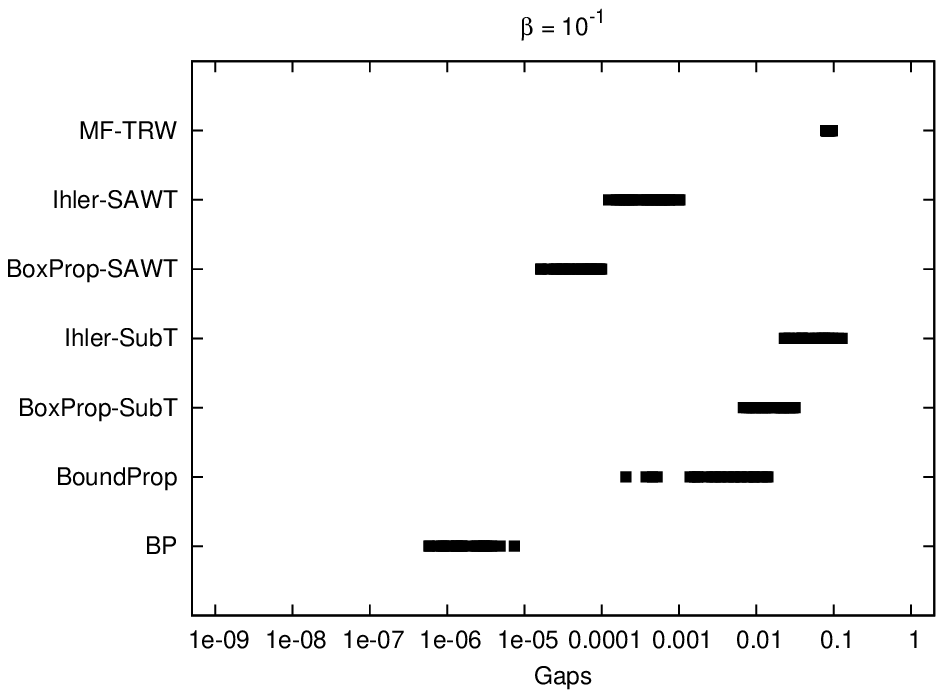}\\[10pt]
\includegraphics[width=0.48\textwidth]{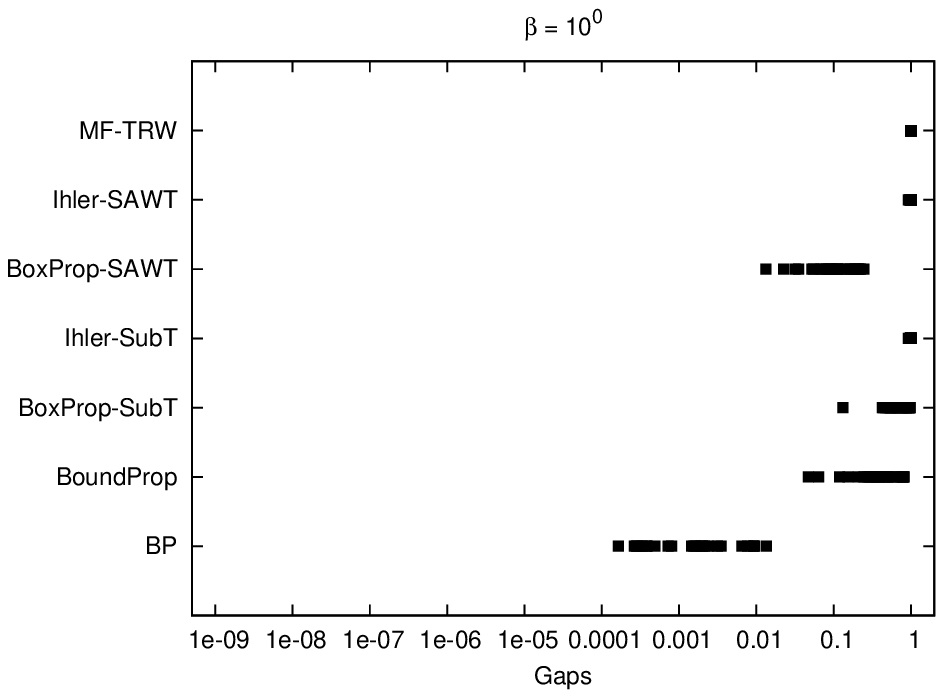}
\includegraphics[width=0.48\textwidth]{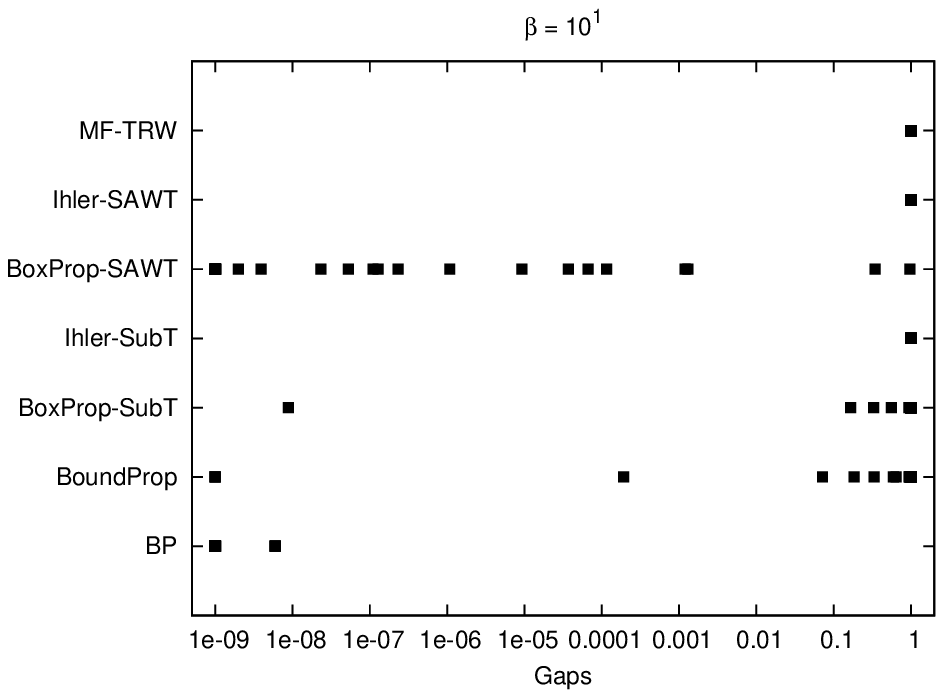}\\[10pt]
\includegraphics[scale=0.7]{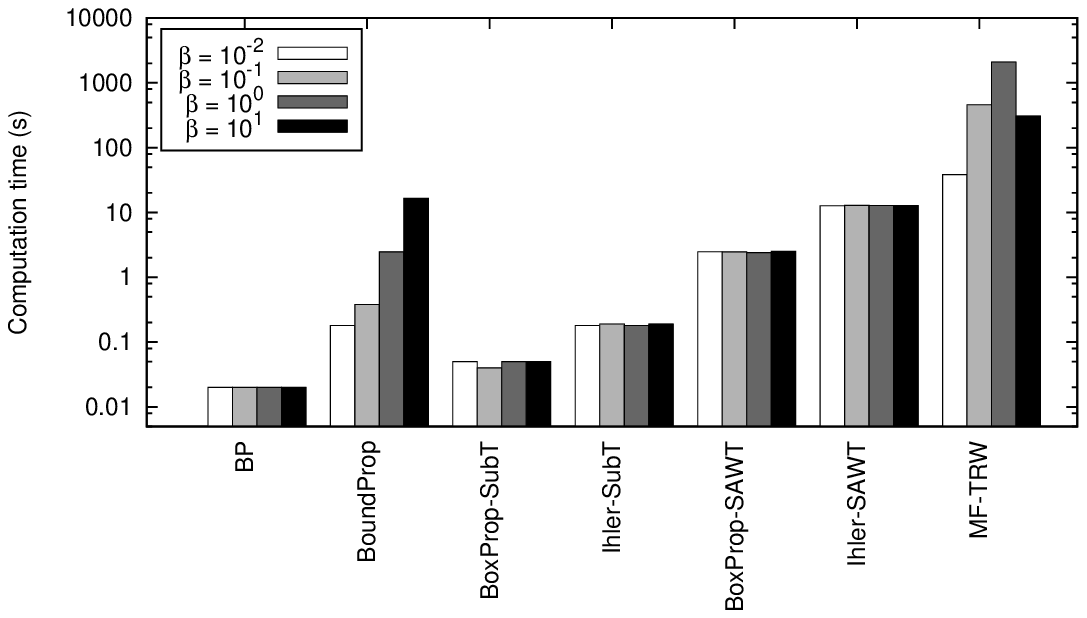}
\caption{Results for grids with ternary variables. The first four graphs show, for different values of 
the interaction strength $\beta$, the gaps of bounds on 
marginals calculated using different methods. Also shown for reference are the errors in the BP 
approximations to the same marginals. The final graph shows the total computation time for each
method.\label{fig:bounds_new_ternary}}
\end{figure}

The results are shown in Figure \ref{fig:bounds_new_ternary}. We have not
compared with bounds involving \JJ\ or \LK\ because these methods have only
been formulated originally for the case of binary variables.  This time, our
method \BoxPropSAW\ yielded the tightest bounds for all interaction strengths
and for all variables (although this is not immediately clear from the plots).

\subsection{Medical diagnosis}

\begin{figure}
\centering
\includegraphics[width=0.45\textwidth]{potent009}\hfill
\includegraphics[width=0.45\textwidth]{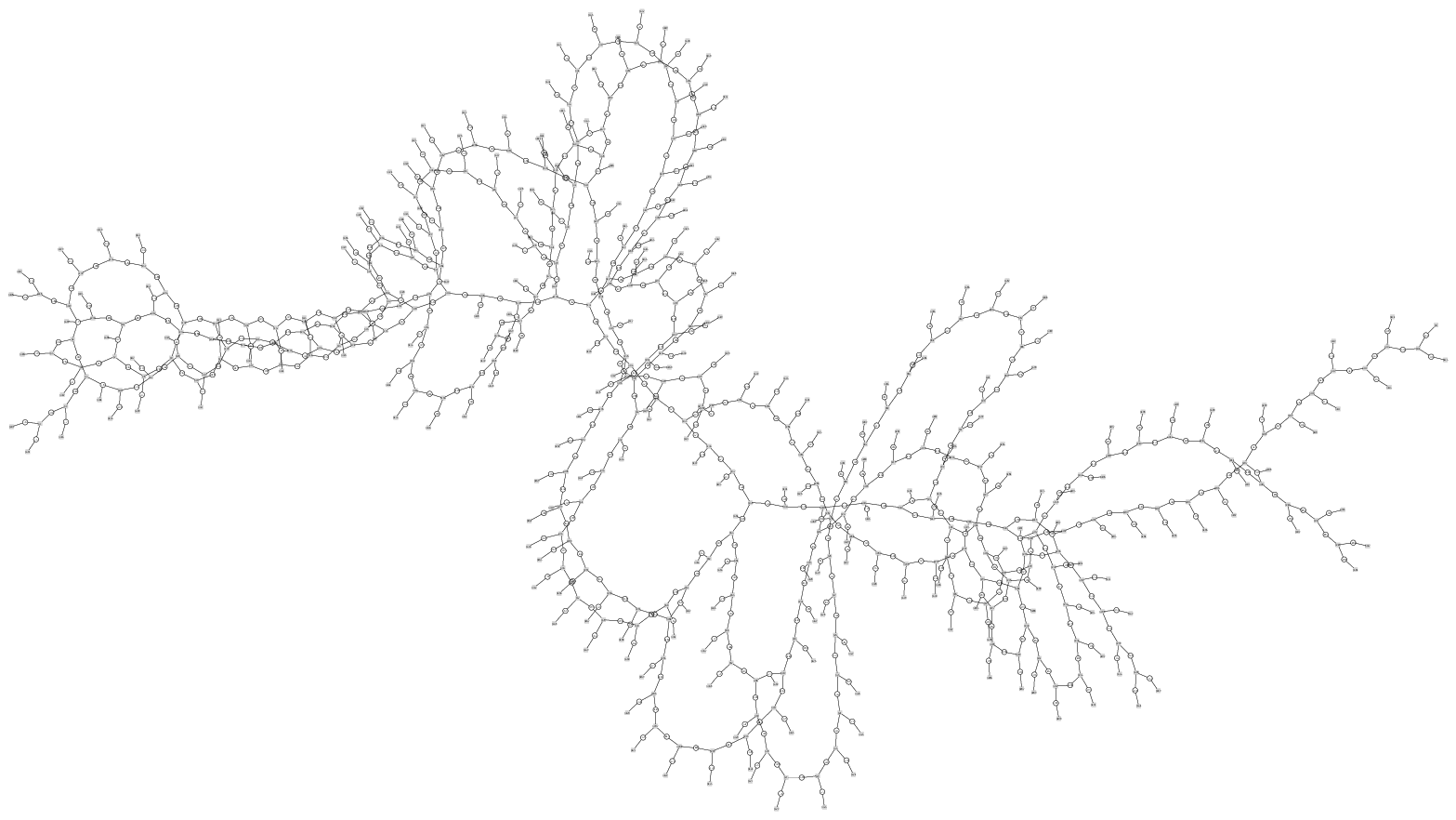}
\caption{Two of the {\sc Promedas} factor graphs.\label{fig:promedas_factor_graphs}}
\end{figure}

We also applied the bounds on simulated {\sc Promedas} patient cases
\citep{WemmenhoveMooijWiegerinckLeisinkKappenNeijt_AIME_07}. These factor
graphs have binary variables and singleton, pairwise and triple interactions
(containing zeros). Two examples of such factor graphs are shown in Figure
\ref{fig:promedas_factor_graphs}. Because of the triple interactions, less methods
were available for comparison.

\begin{figure}
\centering
\includegraphics[width=0.48\textwidth]{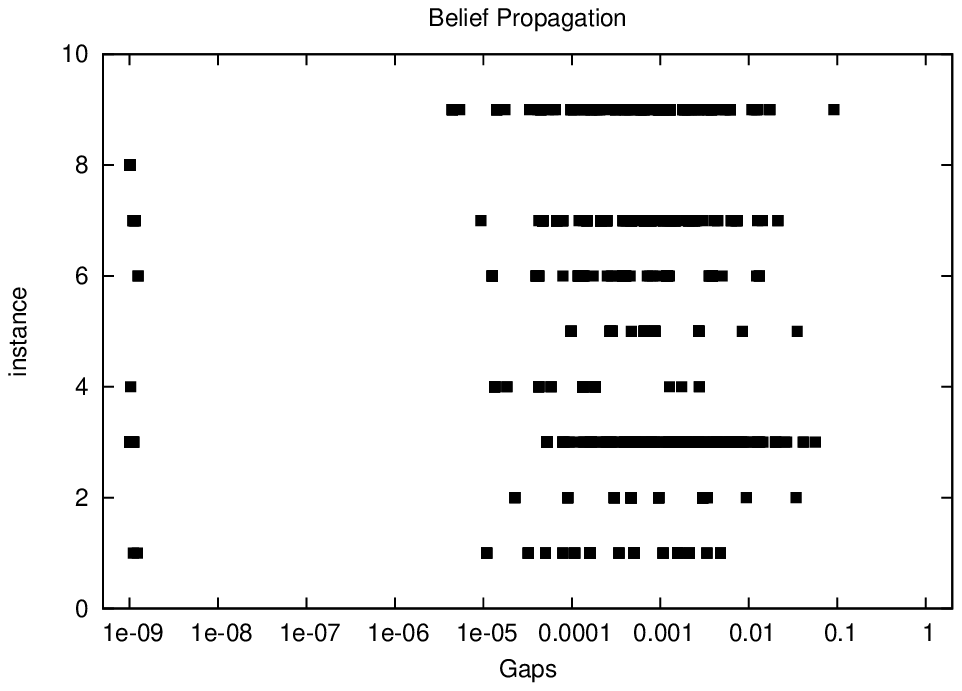}
\includegraphics[width=0.48\textwidth]{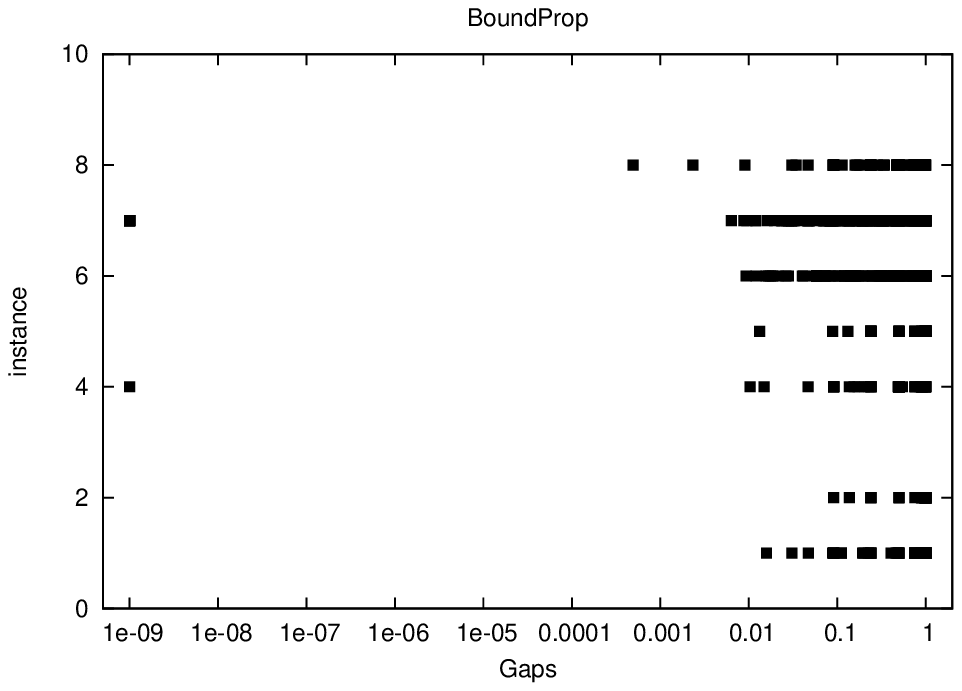}\\[10pt]
\includegraphics[width=0.48\textwidth]{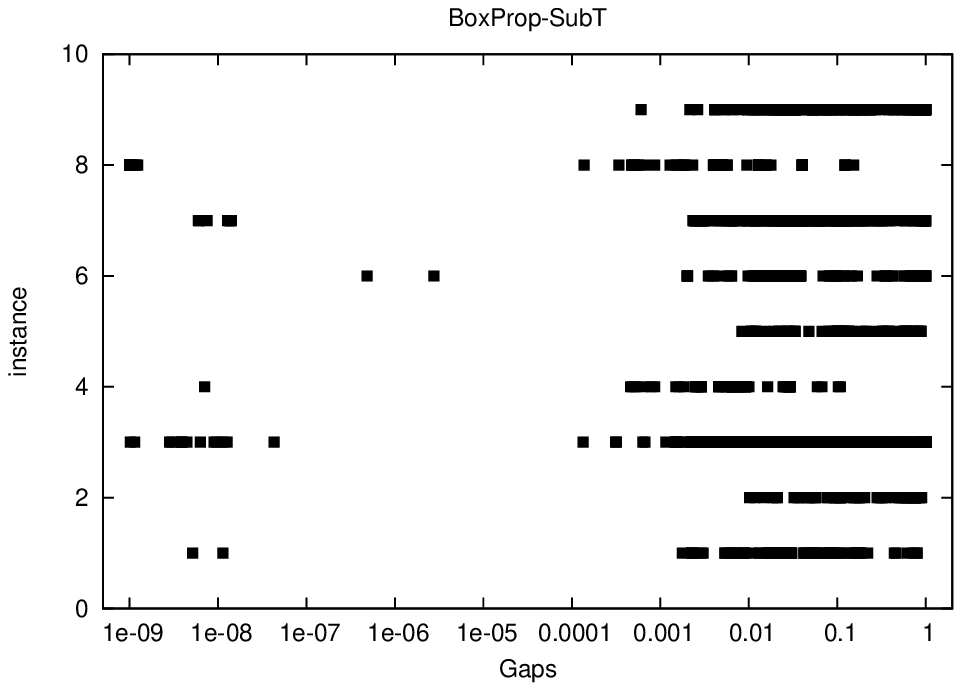}
\includegraphics[width=0.48\textwidth]{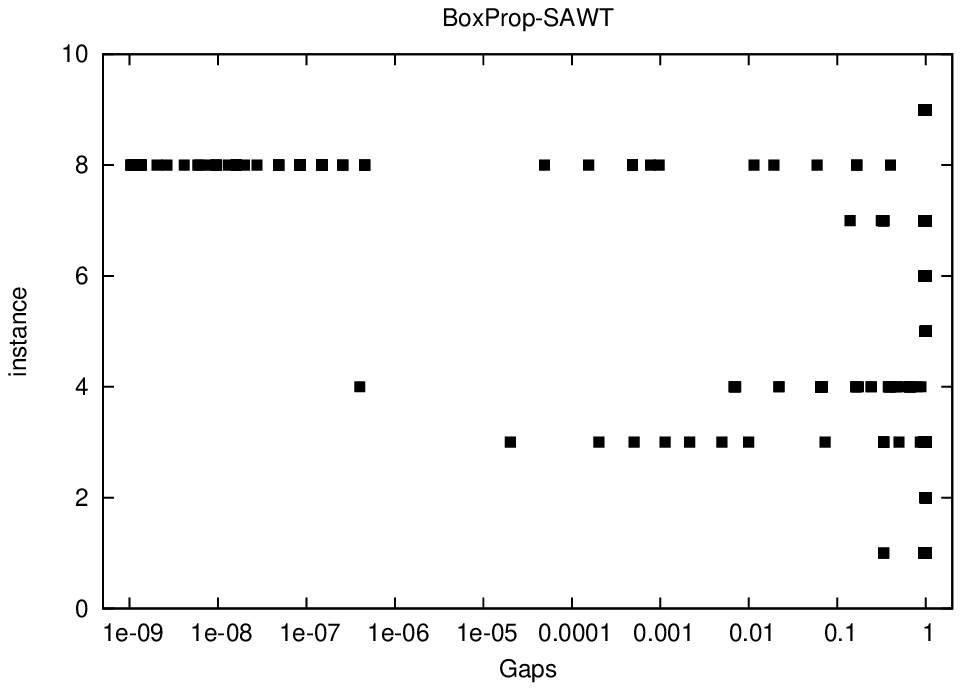}
\caption{Results for nine different factor graphs corresponding
to simulated {\sc Promedas} patient cases. In reading order: Belief
Propagation errors, \BoundProp, \BoxPropSubtree\ and \BoxPropSAW.\label{fig:promedas_bounds}}
\end{figure}

The results of various bounds for nine different, randomly generated,
instances are shown in Figure \ref{fig:promedas_bounds}.
The total number of variables for these nine instances was 1270. The total
computation time needed for \BoxPropSubtree\ was $51\,\mathrm{s}$, for \BoxPropSAW\ $149\,\mathrm{s}$,
for \BoundProp\ more than $75000\,\mathrm{s}$ (we aborted the method for two instances
because convergence was very slow, which explains the missing results in the
plot) and to calculate the Belief Propagation errors took $254\,\mathrm{s}$. 
\BoundProp\ gave the tightest bound for only 1 out of 1270 variables, 
\BoxPropSAW\ for 5 out of 1270 variables and \BoxPropSubtree\ gave the 
tightest bound for the other 1264 variables.

Interestingly, whereas for pairwise interactions, \BoxPropSAW\ gives tighter bounds 
than \BoxPropSubtree, for the factor graphs considered here, the bounds
calculated by \BoxPropSAW\ were generally less tight than those calculated by
\BoxPropSubtree. This is presumably due to the local bound
\eref{eq:boxprop_saw_factor_to_variable} needed on the SAW tree, which is quite
loose compared with the local bound
\eref{eq:boxprop_subtree_factor_to_variable} that assumes independent incoming
bounds. 

Not only does \BoxPropSubtree\ give the tightest bounds for almost all
variables, it is also the fastest method. 
Finally, note that the tightness of these bounds still varies widely depending
on the instance and on the variable of interest.

\section{Conclusion and discussion}\label{sec5:discussion}

We have derived two related novel bounds on exact single-variable marginals.
Both bounds also bound the approximate Belief Propagation marginals. The bounds
are calculated by propagating convex sets of measures over a subtree of the
computation tree, with update equations resembling those of BP. For variables
with a limited number of possible values, the bounds can be computed
efficiently. Empirically, our bounds often outperform the existing
state-of-the-art in that case. Although we have only shown results for factor
graphs for which exact inference was still tractable (in order to be able to
calculate the BP error), we would like to stress here that it is not difficult
to construct factor graphs for which exact inference is no longer tractable 
but the bounds can still be calculated efficiently. An example are large Ising
grids (of size $m\times m$ with $m$ larger than 30). Indeed, for binary Ising 
grids, the computation time of the bounds (for all variables in the network) 
scales linearly in the number of variables, assuming that we truncate the 
subtrees and SAW trees to a fixed maximum size.

Whereas the results of different approximate inference methods usually cannot
be combined in order to get a better estimate of marginal probabilities, for
bounds one can combine different methods simply by taking the tightest bound or
the intersection of the bounds. Thus it is generally a good thing to have
different bounds with different properties (such as tightness and computation
time).

An advantage of our methods \BoxPropSubtree\ and \BoxPropSAW\ over iterative
methods like \BoundProp\ and \MF-\TRW\ is that the computation time of the
iterative methods is difficult to predict (since it depends on the number of
iterations needed to converge which is generally not known a priori). In
contrast, the computation time needed for our bounds \BoxPropSubtree\ and
\BoxPropSAW\ only depends on the structure of the factor graph (and the chosen
subtree) and is independent of the values of the interactions. Furthermore, by
truncating the tree one can trade some tightness for computation time.

By far the slowest methods turned out to be those combining the upper bound
\TRW\ with a lower bound on the partition sum. The problem here is that \TRW\
usually needs many iterations to converge, especially for stronger
interactions where convergence rate can go down significantly. In order to
prevent exceedingly long computations, we had to hand-tune the convergence
criterion of \TRW\ according to the case at hand.

\BoundProp\ can compete in certain cases with the bounds derived here,
but more often than not it turned out to be rather slow or did not yield very
tight bounds. Although \BoundProp\ also propagates bounding boxes over measures,
it does this in a slightly different way which does not exploit independence
as much as our bounds. On the other hand, it can propagate bounding boxes 
several times, refining the bounds more and more each iteration.

Regarding the related bounds \BoxPropSubtree, \BoxPropSAW\ and \BoxPropSAWI\ we
can draw the following conclusions.  For pairwise interactions and variables
that have not too many possible values, \BoxPropSAW\ is the method of choice,
yielding the tightest bounds without needing too much computation time. The 
bounds are more accurate than the bounds produced by \BoxPropSAWI\
due to the more precise local bound that is used; the difference is largest for
strong interactions. However, the computation time of this more precise local
bound is exponential in the number of possible values of the variables, whereas
the local bound used in \BoxPropSAWI\ is only polynomial in the number of possible
values of the variables.  Therefore, if this number is large, \BoxPropSAW\ may
be no longer applicable in practice, whereas \BoxPropSAWI\ still may be
applicable. If factors are present that depend on more than two variables, it
seems that \BoxPropSubtree\ is the best method to obtain tight bounds,
especially if the interactions are strong.  Note that it is not immediately
obvious how to extend \BoxPropSAWI\ beyond pairwise interactions, so we could
not compare with that method in that case.

This work also raises some new questions and opportunities for future work.
First, the bounds can be used to generalize the improved conditions for
convergence of Belief Propagation that were derived in
\citep{MooijKappen_IEEETIT_07} beyond the special case of binary variables
with pairwise interactions. Second, it may be possible to combine the various
ingredients in \BoundProp, \BoxPropSubtree\ and \BoxPropSAW\ in novel ways
in order to obtain even better bounds. Third, it is an interesting open question
whether the bounds can be extended to continuous variables in some way. Finally,
although our bounds are a step forward in quantifying the error of Belief 
Propagation, the actual error made by BP is often at least one order of magnitude 
lower than the tightness of these bounds. This is due to the fact that (loopy) BP 
cycles information through loops in the factor graph; this cycling apparently
improves the results. The interesting and still unanswered question is why it 
makes sense to cycle information in this way and whether this error reduction
effect can be quantified.


\begin{acks}
The research reported here is part of the Interactive Collaborative
Information Systems (ICIS) project (supported by the Dutch Ministry of Economic
Affairs, grant BSIK03024) and was also sponsored in part by the Dutch
Technology Foundation (STW).  

We thank Wim Wiegerinck for several fruitful 
discussions, Bastian Wemmenhove for providing the {\sc Promedas} test cases, 
and Martijn Leisink for kindly providing his implementation of Bound Propagation.
\end{acks}

\vskip 0.2in

\begin{thebibliography}{16}
\providecommand{\natexlab}[1]{#1}
\providecommand{\url}[1]{\texttt{#1}}
\expandafter\ifx\csname urlstyle\endcsname\relax
  \providecommand{\doi}[1]{doi: #1}\else
  \providecommand{\doi}{doi: \begingroup \urlstyle{rm}\Url}\fi

\bibitem[Cooper(1990)]{Cooper90}
G.F. Cooper.
\newblock The computational complexity of probabistic inferences.
\newblock \emph{Artificial Intelligence}, 42\penalty0 (2-3):\penalty0 393--405,
  March 1990.

\bibitem[Georgii(1988)]{Georgii88}
H.-O. Georgii.
\newblock \emph{{G}ibbs Measures and Phase Transitions}.
\newblock Walter de Gruyter, Berlin, 1988.

\bibitem[Ihler(2007)]{Ihler07}
A.~Ihler.
\newblock Accuracy bounds for belief propagation.
\newblock In \emph{Proceedings of the 23th Annual Conference on Uncertainty in
  Artificial Intelligence (UAI-07)}, July 2007.

\bibitem[Jaakkola and Jordan(1996)]{JaakkolaJordan96}
T.~S. Jaakkola and M.~Jordan.
\newblock Recursive algorithms for approximating probabilities in graphical
  models.
\newblock In \emph{Proc. Conf. Neural Information Processing Systems (NIPS 9)},
  pages 487--493, Denver, CO, 1996.

\bibitem[Kschischang et~al.(2001)Kschischang, Frey, and
  Loeliger]{KschischangFreyLoeliger01}
F.~R. Kschischang, B.~J. Frey, and H.-A. Loeliger.
\newblock Factor graphs and the sum-product algorithm.
\newblock \emph{IEEE Trans. Inform. Theory}, 47\penalty0 (2):\penalty0
  498--519, February 2001.

\bibitem[Leisink and Kappen(2003)]{LeisinkKappen03}
M.~Leisink and B.~Kappen.
\newblock Bound propagation.
\newblock \emph{Journal of Artificial Intelligence Research}, 19:\penalty0
  139--154, 2003.

\bibitem[Leisink and Kappen(2001)]{LeisinkKappen01}
M.~A.~R. Leisink and H.~J. Kappen.
\newblock A tighter bound for graphical models.
\newblock In Lawrence~K. Saul, Yair Weiss, and {L\'{e}on} Bottou, editors,
  \emph{Advances in Neural Information Processing Systems 13 (NIPS*2000)},
  pages 266--272, Cambridge, MA, 2001. MIT Press.

\bibitem[Mooij and Kappen(2007)]{MooijKappen_IEEETIT_07}
J.~M. Mooij and H.~J. Kappen.
\newblock Sufficient conditions for convergence of the sum-product algorithm.
\newblock \emph{IEEE Transactions on Information Theory}, 53\penalty0
  (12):\penalty0 4422--4437, December 2007.
\newblock \doi{10.1109/TIT.2007.909166}.

\bibitem[Pearl(1988)]{Pearl88}
J.~Pearl.
\newblock \emph{Probabilistic Reasoning in Intelligent systems: Networks of
  Plausible Inference}.
\newblock Morgan Kaufmann, San Francisco, CA, 1988.

\bibitem[Scott and Sokal(2005)]{ScottSokal05}
A.~D. Scott and A.~D. Sokal.
\newblock The repulsive lattice gas, the independent-set polynomial, and the
  lovasz local lemma.
\newblock \emph{Journal of Statistical Physics}, 118:\penalty0 1151--1261,
  2005.

\bibitem[Taga and Mase(2006)]{TagaMase06b}
Nobuyuki Taga and Shigeru Mase.
\newblock Error bounds between marginal probabilities and beliefs of loopy
  belief propagation algorithm.
\newblock In \emph{MICAI}, pages 186--196, 2006.
\newblock URL \url{http://dx.doi.org/10.1007/11925231_18}.

\bibitem[Tatikonda(2003)]{Tatikonda03}
S.~C. Tatikonda.
\newblock Convergence of the sum-product algorithm.
\newblock In \emph{Proceedings 2003 IEEE Information Theory Workshop}, pages
  222--225, April 2003.

\bibitem[Tatikonda and Jordan(2002)]{TatikondaJordan02}
S.~C. Tatikonda and M.~I. Jordan.
\newblock Loopy belief propogation and {G}ibbs measures.
\newblock In \emph{Proc. of the 18th Annual Conf. on Uncertainty in Artificial
  Intelligence (UAI-02)}, pages 493--500, San Francisco, CA, 2002. Morgan
  Kaufmann Publishers.

\bibitem[Wainwright et~al.(2005)Wainwright, Jaakkola, and
  Willsky]{WainwrightJaakkolaWillsky05}
M.~J. Wainwright, T.~Jaakkola, and A.~S. Willsky.
\newblock A new class of upper bounds on the log partition function.
\newblock \emph{IEEE Transactions on Information Theory}, 51:\penalty0
  2313--2335, July 2005.

\bibitem[Weitz(2006)]{Weitz06}
D.~Weitz.
\newblock Counting independent sets up to the tree threshold.
\newblock In \emph{Proceedings ACM symposium on Theory of Computing}, page
  140–149. ACM, 2006.

\bibitem[Wemmenhove et~al.(2007)Wemmenhove, Mooij, Wiegerinck, Leisink, Kappen,
  and Neijt]{WemmenhoveMooijWiegerinckLeisinkKappenNeijt_AIME_07}
B.~Wemmenhove, J.~M. Mooij, W.~Wiegerinck, M.~Leisink, H.~J. Kappen, and J.~P.
  Neijt.
\newblock Inference in the {P}romedas medical expert system.
\newblock In \emph{Proceedings of the 11th Conference on Artificial
  Intelligence in Medicine (AIME 2007)}, volume 4594 of \emph{Lecture Notes in
  Computer Science}, pages 456--460. Springer, 2007.
\newblock ISBN 978-3-540-73598-4.
\newblock \doi{10.1007/978-3-540-73599-1_61}.

\end{thebibliography}

\clearpage

\end{document}